\documentclass[10pt,twocolumn,paper=letter]{scrartcl}
  
\usepackage{graphicx}
\usepackage{amsmath}
\usepackage{amsfonts}
\usepackage{amssymb}

\usepackage[numbers,sort]{natbib}
\usepackage{xspace}
\usepackage{xcolor}
\usepackage{amssymb}
\usepackage{enumerate}
\usepackage{verbatim}
\usepackage{subfigure}
\usepackage{bbm}
\usepackage{xr}
\usepackage{threeparttable}

\usepackage[colorlinks,linkcolor=blue,citecolor=red]{hyperref}
\usepackage{hypernat}
\usepackage[bf,format=plain,singlelinecheck=false]{caption}

\usepackage{algorithm}
\usepackage{algorithmic}

\newcommand{\ATElasso}{\widehat{ATE}_{\textnormal{Lasso}}\xspace}

\newcommand\ABedit[1]{\textcolor{blue}{#1}}
\newcommand\JSedit[1]{\textcolor{purple}{#1}}
\newcommand\define{\mathrel{\overset{\makebox[0pt]{\mbox{\normalfont\tiny\sffamily def}}}{=}}}
\newcommand\bbeta{\boldsymbol{\beta}}

\renewcommand\ABedit[1]{\textcolor{black}{#1}}
\renewcommand\JSedit[1]{\textcolor{black}{#1}}
\newcommand\CHZedit[1]{\textcolor{black}{#1}}
\newcommand\abedit[1]{\textcolor{black}{#1}}

\newtheorem{thm}{Theorem}
\newtheorem{lem}{Lemma}
\newtheorem{cor}{Corollary}
\newtheorem{condition}{Condition}

\newtheorem{defn}{Def\,inition}
\newcommand{\theoremfont}{\itshape}

\newcommand{\err}{e}  
\newcommand{\st}{\mbox{ s.t. }}
\newcommand{\bx}{\mathbf{x}}
\newcommand{\bh}{\mathbf{h}}
\newcommand{\bz}{\mathbf{z}}
\newcommand{\bw}{\mathbf{w}}
\newcommand{\bu}{\mathbf{u}}
\newcommand{\ATEsig}{\sigma}

\newcommand{\norm}[1]{\left\Vert #1 \right\Vert}
\newcommand{\shortnorm}[1]{\Vert #1 \Vert}

\newenvironment{proof}{%
	\trivlist\item[\hskip \labelsep{\theoremfont Proof.}]%
}{\endtrivlist}

\DeclareMathOperator*{\argmin}{arg\,min}

\renewcommand{\eqref}[1]{\textup{\bf \ref{#1}}}

\title{Lasso adjustments of treatment effect estimates in randomized experiments}

\author{
	Adam Bloniarz\thanks{Department of Statistics, University of California, Berkeley, CA.}
	\and
	Hanzhong Liu\thanks{Department of Statistics, University of California, Berkeley, CA.}
    \and
    Cun-Hui Zhang\thanks{Department of Statistics and Biostatistics, Rutgers University, Piscataway, NJ}
    \and
    Jasjeet Sekhon\thanks{Department of Political Science, Department of Statistics, University of California, Berkeley, CA}
    \and
    Bin Yu\thanks{Department of Statistics, Department of Electrical Engineering and Computer Science, University of California, Berkeley, CA}
}

\begin{document}

	\maketitle

\begin{abstract}
  \noindent	\footnotesize\textbf{Abstract}

  \noindent	We provide a principled way for investigators to analyze randomized experiments when the number of covariates is large. Investigators often use linear multivariate regression to analyze randomized experiments instead of simply reporting the difference of means between treatment and control groups. Their aim is to reduce the variance of the estimated treatment effect by adjusting for covariates. If there are a large number of covariates
  relative to the number of observations, regression may perform poorly because of overfitting. In such cases, the Lasso may be helpful. We study the resulting Lasso-based treatment effect estimator under the Neyman-Rubin model of randomized experiments. We present theoretical conditions that guarantee that the estimator is more efficient than the simple difference-of-means estimator, and we provide a conservative estimator of the asymptotic
  variance, which can yield tighter confidence intervals than the difference-of-means estimator.  Simulation and data examples show that Lasso-based adjustment can be advantageous even when the number of covariates is less than the number of observations. Specifically, a variant using Lasso for selection and OLS for estimation performs particularly well, and it chooses a smoothing parameter based on combined performance of Lasso and OLS.
  \ \\

  \noindent \textbf{Keywords}

  \noindent Randomized experiment, Neyman-Rubin model, average treatment effect, high-dimensional statistics, Lasso, concentration inequality
\end{abstract}

	\section{Introduction}
	
	Randomized experiments are widely used \ABedit{to measure the efficacy of treatments}. Randomization ensures that treatment assignment is not influenced by any potential confounding factors, both observed and unobserved.  Experiments are particularly useful when there is no rigorous theory of a system's dynamics, and full identification of confounders would be impossible.  This advantage was cast elegantly in mathematical terms in the early 20th century by Jerzy Neyman, who introduced a simple model for randomized experiments, which showed that the difference of average outcomes in the treatment and control groups is statistically unbiased for the Average Treatment Effect (ATE) over the experimental sample~\cite{Splawa-Neyman1990}.

However, no experiment occurs in a vacuum of scientific knowledge. Often, \JSedit{baseline covariate information is collected about individuals in an experiment.} \JSedit{Even when treatment assignment is not related to these covariates}, analyses of experimental outcomes often take them into account with the goal of improving \JSedit{the accuracy of treatment effect estimates}.  In modern randomized experiments, the number of covariates can be very large---sometimes even larger than the number of \ABedit{individuals} in the study. In clinical trials overseen by regulatory bodies like the FDA and MHRA, demographic and genetic information may be recorded about each patient. In applications in the tech industry, where randomization is often called A/B testing, there is often a huge amount of behavioral data collected on each user. However, in this `big data' setting, much of this data may be irrelevant to the outcome being studied \JSedit{or there may be more potential covariates than observations, especially once interactions are taken into account. In these cases, selection of important covariates or some form of regularization is necessary for effective regression adjustment.}

To ground our discussion, we examine a randomized trial of the Pulmonary Artery Catheter (PAC) that was carried out in 65 intensive care units in the UK between 2001 and 2004, called PAC-man \cite{harvey2005assessment}. The PAC is a monitoring device commonly inserted into critically ill patients after admission to intensive care, and it provides a continuous measurement of several indicators of cardiac activity. However, insertion of PAC is an invasive procedure that carries some risk of complications (including death), and it involves significant expenditure both in equipment costs and personnel~\cite{Dalen2001}. Controversy over its use came to a head when an observational study found that PAC had an adverse effect on patient survival and led to increased cost of care \cite{connors1996effectiveness}. This led to several large-scale randomized trials, including PAC-man.

In the PAC-man trial, randomization of treatment was largely successful, and a number of covariates were measured about each patient in the study. If covariate interactions are included, the number of covariates exceeds the number of \ABedit{individuals} in the study; however, few of them are predictive of the patient's outcome. As it turned out, the (pre-treatment) estimated probability of death was imbalanced between the treatment and control groups (p = 0.005, Wilcoxon rank sum test). Because the control group had, on average, a slightly higher risk of death, the unadjusted difference-in-means estimator may overestimate the benefits of receiving a PAC.  Adjustment for this imbalance seems advantageous in this case, since the pre-treatment probability of death is clearly predictive of health outcomes post-treatment.

In this paper, we study regression-based adjustment, \ABedit{using the Lasso to select relevant} covariates. Standard linear regression \ABedit{based on ordinary least squares} suffers from over-fitting if a large number of covariates and interaction terms are included in the model. In such cases, researchers sometimes perform model selection based on observing which covariates are unbalanced given the realized randomization. This generally leads to misleading inferences because of incorrect test levels~\cite{permutt1990testing}. The Lasso~\cite{Tibshirani1994} provides researchers with an alternative that can mitigate these problems and still perform model selection. We define an estimator, $\ATElasso$, which is \ABedit{based on} running an $l_1$-penalized linear regression of the outcome on treatment, covariates and, following the method introduced in~\cite{Lin2013}, treatment $\times$ covariate interactions. Because of the geometry of the $l_1$ penalty, the Lasso will usually set many regression coefficients to 0, and is well defined even if the number of covariates is larger than the number of observations. The Lasso's theoretical properties under the standard linear model have been widely studied in the last decade; consistency properties for coefficient estimation, model selection, and out-of-sample prediction are well understood (see \cite{buhlmann2011statistics} for an overview).

In \ABedit{the theoretical analysis in} this paper, instead of assuming that the standard linear model is the true data-generating mechanism, we work under the aforementioned non-parametric model of randomization introduced by Neyman~\cite{Splawa-Neyman1990} and popularized by Donald Rubin~\cite{Rubin1974}. In this model, the outcomes and covariates are fixed quantities, and the treatment group is assumed to be sampled without replacement from a finite population. The treatment indicator, rather than an error term, is the source of randomness, and it determines which of two potential outcomes is revealed to the experimenter. Unlike the standard linear model, the Neyman-Rubin model makes few assumptions not guaranteed by the randomization itself. The setup of the model does rely on the stable unit treatment value assumption (SUTVA), which states that there is only one version of treatment, and that the potential outcome of one unit should be unaffected by the particular assignment of treatments to the other units; however it makes no assumptions of linearity or exogeneity of error terms. Ordinary Least Squares (OLS) \cite{Freedman2008a}\cite{Freedman2008b}\cite{Lin2013}, logistic regression \cite{Freedman2008c}, and post-stratification \cite{Miratrix2013} are among the adjustment methods that have been studied under this model.

To be useful to practitioners, the Lasso-based treatment effect estimator must be consistent and yield a method to construct valid confidence intervals. We outline conditions on the covariates and potential outcomes that will guarantee these properties. \ABedit{We show that an upper bound} for the asymptotic variance can be estimated from the model residuals, yielding asymptotically conservative confidence intervals for the average treatment effect which can be substantially narrower than the unadjusted confidence intervals.  Simulation studies are provided to show the advantage of the Lasso adjusted estimator and to show situations where it breaks down. We \ABedit{apply the} estimator to the PAC-man data, and compare the estimates and confidence intervals derived from the unadjusted, OLS-adjusted, and Lasso-adjusted methods. We also compare different methods of selecting the Lasso tuning parameter on this data.

\section{Framework and definitions}
We give a brief outline of the Neyman-Rubin model for a randomized experiment; the reader is urged to consult \cite{Splawa-Neyman1990}, \cite{Rubin1974}, and \cite{Holland1986} for more details. We follow the notation introduced in \cite{Freedman2008a} and \cite{Lin2013}. For concreteness, we illustrate the model
in the context of the PAC-man trial.

\ABedit{For each individual in the study, the model assumes that there exists a pair of quantities representing his/her health outcomes under the possibilities of receiving and not receiving the catheter. These are called the potential outcomes under treatment and control, and are denoted as $a_i$ and $b_i$, respectively. In the course of the study, the experimenter observes only one of these quantities for each individual, since the catheter is either inserted or not.  The causal effect of the treatment on individual $i$ is defined, in theory, to be $a_i - b_i$, but this is unobservable. Instead of trying to infer individual-level effects, we will assume that the intention is to estimate the average causal effect over the whole population, as outlined in the next section.}

\ABedit{In the mathematical specification of this model we consider the potential outcomes to be fixed, non-random quantities, even though they are not all observable.} The only randomness in the model comes from the assignment of treatment, which is controlled by the experimenter. We define random treatment indicators $T_i$, which take on a value $1$ for a treated individual, or $0$ for an untreated individual.  We will assume that the set of treated individuals is \ABedit{sampled} without replacement from the full population, where the size of the treatment group is fixed beforehand; thus the $T_i$ are identically distributed but not independent.
The model for the observed outcome for individual $i$, defined as $Y_i$, is thus
\begin{equation*}
Y_i = T_i a_i + (1-T_i) b_i.
\end{equation*}
\ABedit{This equation simply formalizes the idea that the experimenter observes the potential outcome under treatment for those who receive the treatment, and the potential outcome under control for those who do not.}

\ABedit{Note that the model does not incorporate any covariate information about the individuals in the study, such as physiological characteristics or health history. However, we will assume we have measured a vector of baseline, pre-experimental covariates for each individual $i$. These might include, for example, age, gender, and genetic makeup. We denote the covariates for individual $i$ as the column vector $\bx_i = (x_{i1},...,x_{ip})^T \in \mathbb{R}^{p}$ and the full design matrix of the experiment as $X=(\bx_1,...,\bx_n)^T$. In the theoretical results, we will assume that there is a correlational relationship between an individual's potential outcomes and covariates, but we will not assume a generative statistical model.}

Define the set of treated individuals as $A = \{i \in \{1,...,n\}:  T_i = 1\}$, and similarly define
\CHZedit{the set of control individuals as $B$}. Define the number of treated and control individuals as  $n_A=\left|A\right|$ and $n_B=\left|B\right|$, respectively, \CHZedit{so that $n_A+n_B=n$}. To indicate averages of quantities over these individuals, we introduce the notation $\bar{\cdot}_A $ and $\bar{\cdot}_B  $. Thus, for example, the average value of the potential outcomes and the covariates in the treatment group are
\[
\bar{a}_A =  n_A^{-1}\hbox{$\sum_{i \in A}$} a_i, \ \bar{\bx}_A  = n_A^{-1} \hbox{$\sum_{i \in A}$}\bx_i,
\]
respectively.
Note that these are random quantities in this model, since the set $A$ is determined by the random treatment assignment. When we want to take the average over the whole population, we will use the notation $\bar{\cdot}$, such as
\[
\bar{a} = n^{-1} \hbox{$\sum_{i=1}^{n}$} a_i, \ \bar{b}   = n^{-1} \hbox{$\sum_{i=1}^{n}$}b_i, \ \bar{\bx}  = n^{-1} \hbox{$\sum_{i=1}^{n}$}\bx_i.
\]
Note that the averages of potential outcomes over the whole population are not considered random, but are unobservable.



\section{Treatment effect estimation}
Our main inferential goal will be average effect of the treatment over the whole population in the study. \ABedit{In a trial such as PAC-man, this represents the difference between the average outcome if everyone had received the catheter, and the average outcome if no one had received it.} This is defined as
\begin{equation*}
ATE = \bar{a} - \bar{b}.
\end{equation*}
The most
natural estimator arises by replacing the population averages with the sample averages:
\begin{equation*}
\label{unadj}
\widehat{ATE}_\textnormal{unadj} \define \bar{a}_A - \bar{b}_B,
\end{equation*}
The subscript ``unadj" indicates an estimator without regression adjustment. The foundational work in \cite{Splawa-Neyman1990} points out that, under a randomized assignment of treatment, $\widehat{ATE}_\textnormal{unadj}$ is unbiased for $ATE$, and derives a conservative procedure for estimating its variance.

While $\widehat{ATE}_\textnormal{unadj}$ is an attractive estimator, covariate information can be used to make adjustments in the hope of reducing variance. A commonly used estimator is
\begin{align*}
\widehat{ATE}_\textnormal{adj} =
\left[ \bar{a}_A - \left( \bar{\bx}_A - \bar{\bx} \right)^T \hat \bbeta^{(a)} \right]  - \left[ \bar{b}_B - \left( \bar{\bx}_B - \bar{\bx} \right)^T \hat \bbeta^{(b)} \right]
\end{align*}
where $ \hat \bbeta^{(a)}, \hat \bbeta^{(b)} \in \mathbb{R}^p$ are adjustment vectors for the treatment and control groups, \ABedit{respectively, as indicated by the superscripts}. The terms $ \bar{\bx}_A - \bar{\bx} $ and $ \bar{\bx}_B - \bar{\bx} $ represent the fluctuation of the covariates in the subsample relative to the full sample, and the adjustment vectors
fit 
the linear relationships between the covariates and
potential outcomes under treatment and control. \ABedit{For example, in the PAC-man trial, this would help alleviate the imbalance in the pre-treatment estimated probability of death: the corresponding element of $\bar{\bx}_B - \bar{\bx}$ would be positive (due to the higher average probability of death in the control group), the corresponding element of $\hat \bbeta^{(b)}$ would be negative (a higher probability of death correlates with worse health outcomes), so the overall treatment effect estimate would be adjusted downwards.}
 This procedure 
is equivalent to imputing the unobserved potential outcomes; if we define
\begin{align*}
\hat{ \bar{a}}_B = \bar{a}_A + \left( \bar{\bx}_B - \bar{\bx}_A \right)^T \hat \bbeta^{(a)},\ 
\hat{ \bar{b}}_A = \bar{b}_B + \left( \bar{\bx}_A - \bar{\bx}_B \right)^T \hat \bbeta^{(b)},
\end{align*}
we can form the equivalent estimator
\begin{equation*}
\widehat{ATE}_\textnormal{adj} = n^{-1}\left( n_A \bar{a}_A + n_B \hat{ \bar{a}}_B \right) -
n^{-1}\left( n_B \bar{b}_B + n_A \hat{ \bar{b}}_A \right).
\end{equation*}
If we consider these adjustment vectors to be fixed (non-random), \ABedit{or if they are derived from an independent data source}, then this estimator is still unbiased, and may have substantially smaller asymptotic and finite-sample variance than the unadjusted estimator. This allows for construction of tighter confidence intervals for the true treatment effect.

In practice, the ``ideal" linear adjustment vectors, \ABedit{leading to a minimum-variance estimator of the form of $\widehat{ATE}_\textnormal{adj}$}, \ABedit{cannot be computed} from the observed data. However, they can be estimated, possibly at the expense of introducing modest finite-sample bias into the treatment effect estimate.  In the classical setup, when the number of covariates is relatively small, ordinary least squares (OLS) regression can be used. The asymptotic properties of this kind of estimator are explored under the Neyman-Rubin model in \cite{Freedman2008b}, \cite{Freedman2008c}, and \cite{Lin2013}. We will follow a particular scheme which is studied in \cite{Lin2013} and shown to have favorable properties: we regress the outcome on treatment indicators, covariates, and treatment $\times$ covariate interactions. This is equivalent to running separate regressions in the treatment and control groups of outcome against an intercept and covariates. If we define $\hat\bbeta^{(a)}_{\textnormal{OLS}}$ and $\hat\bbeta^{(b)}_{\textnormal{OLS}}$ as the coefficients from the separate regressions, then the estimator is
\begin{align*}
\widehat{ATE}_\textnormal{OLS}  = &
\left[ \bar{a}_A - \left( \bar{\bx}_A - \bar{\bx} \right)^T\hat\bbeta^{(a)}_{\textnormal{OLS}} \right] \\
& - \left[ \bar{b}_B - \left( \bar{\bx}_B - \bar{\bx} \right)^T\hat\bbeta^{(b)}_{\textnormal{OLS}} \right].
\end{align*}
This has some finite-sample bias, but \cite{Lin2013} shows that it vanishes quickly at the rate of $1/n$ under moment conditions on the potential outcomes
and covariates. Moreover, for a fixed $p$, under regularity conditions, the inclusion of interaction terms guarantees that it never has higher asymptotic variance than the unadjusted estimator, and asymptotically conservative confidence intervals for the true parameter can be constructed.

\ABedit{In modern randomized trials, where a large number of covariates are recorded for each individual, $p$ may be comparable to or even larger than $n$. In this case OLS regression can overfit the data badly, or may even be ill-posed, leading to estimators with large finite-sample variance.} To remedy this, we propose estimating the adjustment vectors using the Lasso \cite{Tibshirani1994}. The adjustment vectors would take the form
{\small \begin{equation}
\begin{split}
\label{def-lasso-a}
\hat{\bbeta}_{\textnormal{Lasso}}^{(a)} = \argmin_{\bbeta} \biggl[
\frac{1}{2n_A} &\sum_{i \in A} \left( a_i - \bar{a}_A - (\bx_i -  \bar{\bx}_A )^T \bbeta \right)^2 \biggr.  \\
& \biggl. + \lambda_a \sum_{j=1}^{p} |\beta_j| \biggr],
\end{split}
\end{equation}}
{\small \begin{equation}
\begin{split}
\label{def-lasso-b}
\hat{\bbeta}_{\textnormal{Lasso}}^{(b)} = \argmin_{\bbeta} \biggl[
\frac{1}{2n_B} &\sum_{i \in B} \left( b_i - \bar{b}_B - (\bx_i -  \bar{\bx}_B )^T \bbeta \right)^2 \biggr.
\\ & \biggl. + \lambda_b \sum_{j=1}^{p} |\beta_j| \biggr],
\end{split}
\end{equation}}
and the proposed Lasso adjusted ATE estimator is\footnote{To simplify the notation, we omit the dependence of $\hat{\bbeta}_{\textnormal{Lasso}}^{(a)}$, $\hat{\bbeta}_{\textnormal{Lasso}}^{(b)}$, $\lambda_a$ and $\lambda_b$ on the population size $n$.}\label{ft:depend-on-n}
\begin{equation}
\begin{split}
\label{def-ate-lasso}
\widehat{ATE}_\textnormal{Lasso} = &
\left[ \bar{a}_A - \left( \bar{\bx}_A - \bar{\bx} \right)^T\hat\bbeta^{(a)}_{\textnormal{Lasso}} \right]
\\ \quad
& - \left[ \bar{b}_B - \left( \bar{\bx}_B - \bar{\bx} \right)^T\hat\bbeta^{(b)}_{\textnormal{Lasso}} \right].  \nonumber
\end{split} \nonumber
\end{equation}
\CHZedit{Here $\lambda_a$ and $\lambda_b$ are regularization parameters} for the Lasso which must be chosen by the experimenter; simulations show that cross-validation works well. In the next section, we study this estimator under the Neyman-Rubin model, and provide conditions on the potential outcomes, the covariates and the regularization parameters under which $\widehat{ATE}_\textnormal{Lasso} $ enjoys similar asymptotic and finite-sample advantages as $\widehat{ATE}_{\textnormal{OLS}}$.

It is worth noting that when two different adjustments are made for the treatment and control groups as in \cite{Lin2013} and here, the covariates
do not have to be the same for the two groups. However, when they are not the same, the Lasso or OLS adjusted estimators are no longer guaranteed
to have smaller or equal asymptotic variance than the unadjusted one, even in the case of fixed $p$.
In practice, one may still choose between the adjusted and unadjusted estimators based on
the widths of the corresponding confidence intervals.



\section{Theoretical results}\label{sec:theory}

\subsection{Notation}
For a vector $\bbeta \in R^p$ and a subset $S \subset \{1,...,p\}$, let  $\beta_j$ be the $j$-th component of $\bbeta$, $\bbeta_S = (\beta_j: j \in S)^T$,
$S^c$ be the complement of $S$, and $|S|$ the cardinality of the set $S$.
For any column vector $\mathbf{u}=(u_1,...,u_m)^T$, let $\|\mathbf{u}\|_2^2 = \sum_{i=1}^{m} u_i^2$, $\|\mathbf{u}\|_1 = \sum_{i=1}^{m} |u_i|$, $\|\mathbf{u}\|_\infty = \max_{i=1,\ldots,m} |u_i|$ and $\|\mathbf{u}\|_0 = | \{j: u_j \neq 0 \} |$. For a given $m\times m$ matrix $D$, let $\lambda_{\textnormal{min}} (D)$ and $\lambda_{\textnormal{max}} (D)$ be the smallest and largest eigenvalues of $D$
respectively, and $D^{-1}$ the inverse of the matrix $D$. Let $\stackrel{d}{\rightarrow}$ and $\stackrel{p}{\rightarrow}$ denote convergence in distribution and in probability, respectively.

\subsection{Decomposition of the potential outcomes}\

The Neyman-Rubin model does not assume a linear relationship
  between the potential outcomes and the covariates. In order to study
  the properties of adjustment under this model, we decompose the
  potential outcomes into a term linear in the covariates and an
  error term. Given vectors of coefficients $\bbeta^{(a)},
\bbeta^{(b)} \in \mathbb{R}^p$, we write\footnote{Again, we omit the
  dependence of $\bbeta^{(a)}$, $\bbeta^{(b)}$, $\lambda_a$,
  $\lambda_b$, $\err^{(a)}$ and $\err^{(b)}$ on $n$.} for $i=1,...,n$,
\begin{equation}
\label{decom-a}
a_i = \bar{a} + ( \bx_i - \bar{\bx} )^T \bbeta^{(a)} + \err^{(a)}_i ,
\end{equation}
\begin{equation}
\label{decom-b}
b_i = \bar{b} + ( \bx_i - \bar{\bx} )^T \bbeta^{(b)} + \err^{(b)}_i.
\end{equation}

Note that we have not added any assumptions to the model; we have simply defined unit-level residuals, $\err_i^{(a)}$ and  $\err_i^{(b)}$, given the vectors $\bbeta^{(a)}, \bbeta^{(b)}$. All the quantities in \eqref{decom-a} and \eqref{decom-b} are fixed, deterministic numbers.
It is easy to verify that $\bar{e}^{(a)} =  \bar{e}^{(b)} =0$.
In order to pursue a theory for the Lasso, we will add assumptions on the populations of $a_i$'s, $b_i$'s, and $\bx_i$'s, and we will assume the existence of $\bbeta^{(a)}, \bbeta^{(b)}$ such that the error terms satisfy certain assumptions.

\subsection{Conditions}
We will need the following to hold for both the treatment and control potential outcomes. The first set of assumptions (\ref{first-cond}-\ref{cond:limit}) are similar to those found in \cite{Lin2013}.
\begin{condition}\label{first-cond}
Stability of treatment assignment probability. 
\begin{eqnarray}
\label{cond:treat-prob}
n_A/n 
\rightarrow p_A, \ \textnormal{as} \ n \rightarrow \infty
\end{eqnarray}
for some $p_A \in \left(0,1\right)$.
\end{condition}

\begin{condition} \label{cond:moment}
The centered moment conditions.
There exists a fixed constant $L>0$ such that, for all $n=1,2,...$ and $j=1,...,p$,
\begin{equation}
\label{cond:xmoment}
 \hbox{$n^{-1} \sum_{i=1}^{n}$} \left( x_{ij}-(\bar{\bx})_j \right)^4 \leq L;
\end{equation}
\begin{equation}
\label{cond:errmoment}
  \hbox{$n^{-1} \sum_{i=1}^{n}$} (e_i^{(a)})^4 \leq L; \ \  \hbox{$n^{-1} \sum_{i=1}^{n}$} (e_i^{(b)})^4 \leq L.
\end{equation}
\end{condition}

\begin{condition} \label{cond:limit}
The means $\hbox{$n^{-1} \sum_{i=1}^{n}$} (e_i^{(a)})^2$,  $\hbox{$n^{-1} \sum_{i=1}^{n}$} (e_i^{(b)})^2$ and $\hbox{$n^{-1} \sum_{i=1}^{n}$} e_i^{(a)} e_i^{(b)}$ converge to finite limits.
\end{condition}


Since we consider the high-dimensional setting where $p$ is
allowed to be 
much larger than $n$, we need additional assumptions to ensure that the Lasso is consistent for estimating $\bbeta^{(a)}$ and $\bbeta^{(b)}$. Before stating them, we define several quantities.

\begin{defn}
Given $\bbeta^{(a)}$ and $\bbeta^{(b)}$, the sparsity measures for treatment and control groups, $s^{(a)}$ and $s^{(b)}$, are defined as the number of nonzero elements of $\bbeta^{(a)}$ and $\bbeta^{(b)}$, i.e.,
\begin{equation}\label{def:s}
s^{(a)} = | \{j: \beta_j^{(a)} \neq 0  \} |, \ s^{(b)} =  | \{j: \beta_j^{(b)} \neq 0  \} |,
\end{equation}
respectively. We will allow $s^{(a)}$ and $s^{(b)}$  to grow with $n$, though the notation does not explicitly show this.
\end{defn}

\begin{defn}
Define $\delta_n$ to be the maximum covariance between the error terms and the covariates.
{\small \begin{equation}\label{def:delta}
\delta_n = \max_{\omega=a,b} \left\{ \max_j  \left|\frac{1}{n} \sum_{i=1}^{n} \left( x_{ij}-(\bar{\bx})_j \right) \left( e^{(\omega)}_i - \bar{e}^{(\omega)} \right) \right| \right\}.
\end{equation}}
\end{defn}

The following conditions will guarantee that the Lasso consistently estimates the adjustment vectors $\bbeta^{(a)}, \bbeta^{(b)}$ at a fast enough rate to ensure asymptotic normality of  $\ATElasso$. It is an open question whether a weaker form of consistency would be sufficient for our results to hold.

\begin{condition} \label{cond:scaling}
Decay and scaling.
Let $s = \max \left\{ s^{(a)}, s^{(b)} \right\}$.
\begin{equation}
\label{cond:delta_n}
\delta_n = o\left( \frac{1}{s\sqrt{\log p}} \right).
\end{equation}
\begin{equation}
\label{cond:s-scaling}
(s \log p)/{\sqrt n} = o(1).
\end{equation}
\end{condition}

\begin{condition}\label{last-cond}
Cone invertibility factor. Define the Gram matrix as
$\Sigma = n^{-1}\sum_{i=1}^n ( \bx_i - \bar{\bx} )( \bx_i - \bar{\bx} )^T$:
There exist constants $C >0$ and $\xi >1$ not depending on $n$, such that
\begin{equation}\label{cond:cone-inv}
\|\bh_S\|_1 \leq Cs \| \Sigma \bh \|_\infty, \ \forall \bh \in \mathcal{C},
\end{equation}
with $\mathcal{C}=\{\bh: \|\bh_{S^c}\|_1 \leq \xi \|\bh_{S}\|_1\}$, and
\begin{equation} \label{def:bigS}
S = \{ j: \bbeta^{(a)}_j \neq 0 \ \textnormal{or} \  \bbeta^{(b)}_j \neq 0 \}.
\end{equation}

\end{condition}

\begin{condition}\label{last-last-cond}
Let $\tau = \min\big\{1/70,(3p_A)^2/70,(3-3p_A)^2/70 \big\}$. For constants $0< \eta < \frac{\xi - 1}{\xi + 1}$ and $\frac{1}{\eta}<M<\infty$, assume the regularization parameters of the Lasso belong to the sets
\begin{equation}\label{cond:lambda-a}
\lambda_a \in   (\frac{1}{\eta}, M] \times \left( \frac{2(1+\tau)L^{1/2}}{p_A} \sqrt{ \frac{2\log p}{n} } + \delta_n \right),
\end{equation}
\begin{equation} \label{cond:lambda-b}
\lambda_b \in   (\frac{1}{\eta}, M] \times \left( \frac{2(1+\tau)L^{1/2}}{p_B} \sqrt{ \frac{2\log p}{n} } + \delta_n \right).
\end{equation}
\end{condition}


Denote respectively the population variances of $e^{(a)}$ and $e^{(b)}$ and the population covariance between them by
\[ \sigma^2_{e^{(a)}} = \hbox{$n^{-1} \sum_{i=1}^{n}$} (e_i^{(a)})^2, \ \
\sigma^2_{e^{(b)}}  = \hbox{$n^{-1} \sum_{i=1}^{n}$} (e_i^{(b)})^2, \]
\[ \sigma_{e^{(a)}e^{(b)}} =   \hbox{$n^{-1} \sum_{i=1}^{n}$} e_i^{(a)} e_i^{(b)}. \]

\begin{thm} \label{normality-thm}
Assume conditions \ref{first-cond} through \ref{last-last-cond} hold for some $ \bbeta^{(a)} $ and $ \bbeta^{(b)}$. Then
\begin{eqnarray} \label{ATE-conv}
\sqrt{n}  \left( \ATElasso - ATE \right) \stackrel{d}{\rightarrow} \mathcal{N}\left( 0, \ATEsig^2 \right)
\end{eqnarray}
where
\begin{equation} \label{sig-def}
\ATEsig^2 = \lim_{n\rightarrow \infty}\left[\frac{1-p_A}{p_A} \sigma^2_{e^{(a)}} + \frac{p_A}{1-p_A}\sigma^2_{e^{(b)}} + 2\sigma_{e^{(a)}e^{(b)}}\right].
\end{equation}
\end{thm}


The proof of Theorem \ref{normality-thm} is given in the supporting information. It is easy to show, as in the following corollary of Theorem \ref{normality-thm}, that the asymptotic variance of $\widehat{ATE}_\textnormal{Lasso}$ is no worse than $\widehat{ATE}_\textnormal{unadj}$ when $\bbeta^{(a)}$ and $\bbeta^{(b)}$ are defined as coefficients of regressing potential outcomes on a subset of covariates. More specifically, suppose there exists a subset $J \subset  \{1,...,p\}$, such that
\begin{equation} \label{beta-a-b}
\bbeta^{(a)} = ( (\bbeta_J^{(a)})^T, \mathbf{0} )^T, \ \bbeta^{(b)} = ( (\bbeta_J^{(b)})^T, \mathbf{0} )^T,
\end{equation}
where $\bbeta_J^{(a)}$ and $\bbeta_J^{(b)}$ are the population level OLS coefficients for regressing the potential outcomes $a$ and $b$ on the covariates in the subset $J$ with intercept, respectively.

\begin{cor}
\label{compare-asym-var}
For $\bbeta^{(a)}$ and $\bbeta^{(b)}$ defined in \eqref{beta-a-b} and some $\lambda_a$ and $\lambda_b$, assume conditions \ref{first-cond} through \ref{last-last-cond} hold.
Then the asymptotic variance of $\sqrt{n} \
\widehat{ATE}_\textnormal{Lasso}$ is no greater than that of the $\sqrt{n} \ \widehat{ATE}_\textnormal{unadj}$. The difference is $\frac{1}{p_A(1-p_A)}\Delta$, where
\begin{eqnarray}
\label{var-compare}
\Delta  =  - \lim_{n\rightarrow \infty} \| X\bbeta_E \|_2^2 \leq 0,
\end{eqnarray}
\begin{eqnarray}
\label{var-compare2}
\bbeta_E =  (1-p_A)\bbeta^{(a)} + p_A \bbeta^{(b)}.
\end{eqnarray}
\end{cor}

\begin{figure}[!ht]
\centerline{\includegraphics[width=.45\textwidth]{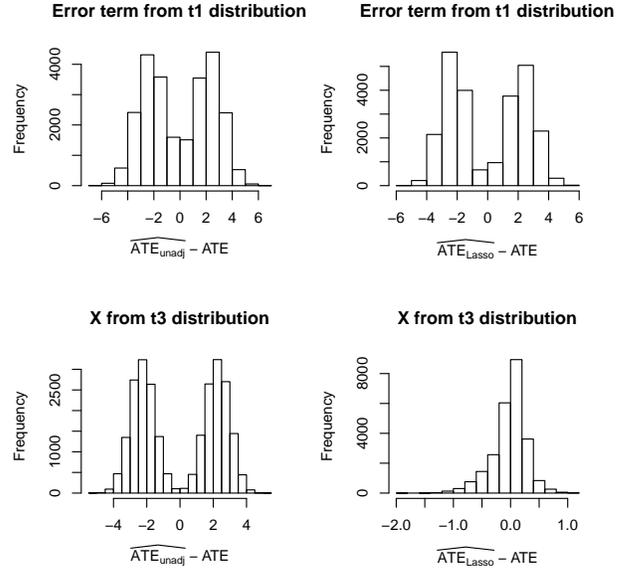}}
\vspace*{0.05in}
\caption{Histograms of the unadjusted estimator and the Lasso adjusted estimator when the moment conditions do not hold. We select the tuning parameters for Lasso using 10-fold cross validation. The potential outcomes are simulated from linear regression model and then kept fixed, see more details in the supporting information. For the upper two subplots, the error terms are generated from $t$ distribution with one degree of freedom and therefore do not satisfy second moment condition; while for the lower two subplots, the covariates are generated from $t$ distribution with there degrees of freedom and thus violate fourth moment condtion.}\label{fig:cond-not-hold}
\end{figure}

\noindent {\bf Remark 1.} {\em If, \abedit{instead of Condition \ref{cond:xmoment}, we} assume that the
  covariates are uniformly bounded, i.e., $\max_{i,j}|x_{ij}|\leq L$,
  then the fourth moment condition on the error terms, given in
  $\eqref{cond:errmoment}$, can be weakened to a second moment
  condition. \abedit{While we do not prove the necessity of any of our conditions, o}ur simulation studies show that the distributions
    of the unadjusted and the Lasso adjusted estimator may be
    non-normal when: (1) The covariates are generated from Gaussian
    distributions and the error terms do not satisfy second moment
    condition, e.g., being generated from a $t$ distribution with one
    degree of freedom; \abedit{or} (2) The covariates do not have bounded
    fourth moments, e.g., being generated from a $t$ distribution with
    three degrees of freedom. See the histograms in
    Figure~\ref{fig:cond-not-hold} where the corresponding p-values of
    Kolmogorov--Smirnov testing for normality are less than
    $2.2e-16$. These findings indicate that our moment conditions
    cannot be dramatically weakened for asymptotic normality. However,
    we also find that the Lasso adjusted estimator still has smaller
    variance and mean squared error than the unadjusted estimator,
    even when these moment conditions do not hold. In practice, when
    the covariates do not have bounded fourth moments, one may perform some transformation---e.g., a logarithm
    transformation---to ensure that the transformed covariates have
    bounded fourth moments while having a sufficiently large variance so as to
    retain useful information. We leave \abedit{it as} future work to explore the
    properties of different transformations.}

\noindent
{\bf Remark 2.} {\em
Statement $\eqref{cond:s-scaling}$, typically required in de-biasing the Lasso \cite{Zhang2014},
is stronger by a factor of $\sqrt{\log p}$
than the usual requirement for $l_1$ consistency of the Lasso.

}

\noindent
{\bf Remark 3.} {\em Condition \ref{last-cond} is slightly weaker than the typical restricted eigenvalue condition 
for analyzing the Lasso.}

\noindent
{\bf Remark 4.} {\em If we assume $\delta_n = O\left(\sqrt{\frac{\log p}{n}}\right)$ which satisfies $\eqref{cond:delta_n}$, then Condition \ref{last-last-cond} requires that the tuning parameters are proportional to $\sqrt{\frac{\log p}{n}}$ which is 
typically assumed for the Lasso in the high-dimensional linear regression model.}

\noindent
{\bf Remark 5.} {\em
For fixed $p$, $\delta_n=0$ in $\eqref{def:delta}$, Condition \ref{cond:scaling} holds automatically,
and Condition \ref{last-cond} holds when the smallest eigenvalue of $\Sigma$ is uniformly bounded
away from 0. In this case,
Corollary \ref{compare-asym-var} reverts to Corollary 1.1. in
\cite{Lin2013}. When these conditions are not satisfied, we should set
$\lambda_a$ and $\lambda_b$ to be large enough to cause the Lasso
adjusted estimator to revert to the unadjusted one.}




\section{Neyman-type conservative variance estimate}

We note that the asymptotic variance in Theorem~\ref{normality-thm} involves the cross-product term $\sigma_{e^{(a)}e^{(b)}}$ which is not consistently estimable in the Neyman-Rubin model as $a_i$ and $b_i$ are never simultaneously observed. However, we can give a Neyman-type conservative estimate of the variance. Let
\begin{equation}
\label{var_estim_a}
 \hat \sigma^2_{e^{(a)}} = \frac{1}{n_A - df^{(a)} } \sum_{i \in A} \left( a_i - \bar{a}_A - (\bx_i - \bar{\bx}_A )^T \hat \bbeta^{(a)}_\textnormal{Lasso} \right)^2,
\end{equation}
\begin{equation}
\label{var_estim_b}
 \hat \sigma^2_{e^{(b)}} = \frac{1}{n_B - df^{(b)} } \sum_{i \in B} \left( b_i - \bar{b}_B - (\bx_i - \bar{\bx}_B )^T \hat \bbeta^{(b)}_\textnormal{Lasso}  \right)^2,
\end{equation}
where $df^{(a)}$ and $df^{(b)}$ are degrees of freedom defined by
\[ df^{(a)} = \hat s^{(a)}+1 = || \hat \bbeta^{(a)}_\textnormal{Lasso} ||_0 +1; \]
\[ df^{(b)} = \hat s^{(b)}+1 = || \hat \bbeta^{(b)}_\textnormal{Lasso} ||_0 +1. \]

Define the variance estimate of $\sqrt{n}(\widehat{ATE}_\textnormal{Lasso}  - ATE )$ as follows:
\begin{equation} \label{sigma-lasso}
\hat \sigma^2_\textnormal{Lasso} = \frac{n}{n_A}  \hat \sigma^2_{e^{(a)}} + \frac{n}{n_B} \hat \sigma^2_{e^{(b)}}.
\end{equation}

\begin{condition}\label{add-cond}
For the Gram matrix $\Sigma$ defined in Condition \ref{last-cond},
the largest eigenvalue is bounded away from $\infty$, that is, there exists a constant $\Lambda_{max} < \infty$ such that
\begin{equation*}
\lambda_{max}\left(\Sigma\right)\leq \Lambda_{max}.
\end{equation*}
\end{condition}

\begin{thm}
\label{conservative_variance}
Assume conditions in Theorem~\ref{normality-thm} and condition~\ref{add-cond} hold. Then $\hat \sigma^2_\textnormal{Lasso}$ converges in probability to
\begin{equation*}
\frac{1}{p_A} \mathop {\lim}\limits_{n \rightarrow \infty } \sigma^2_{e^{(a)}} + \frac{1}{1-p_A} \mathop {\lim}\limits_{n \rightarrow \infty } \sigma^2_{e^{(b)}},
\end{equation*}
which is greater than or equal to the asymptotic variance of $\sqrt{n}(\widehat{ATE}_\textnormal{Lasso}  - ATE )$. The difference is
\begin{equation*}
\begin{split}
\mathop {\lim}\limits_{n \rightarrow \infty } & \frac{1}{n} \sum_{i=1}^{n}   \biggl[ a_i - b_i - ATE  - (\bx_i - \bar{\bx} )^T (\bbeta^{(a)} - \bbeta^{(b)}) \biggr]^2.
\end{split}
\end{equation*}

\end{thm}

\noindent
{\bf Remark 6.} {\em The Neyman-type conservative variance estimate for the unadjusted estimator is given by
{ \small $$ \hat \sigma^2_\textnormal{unadj} = \frac{n}{n_A} \frac{1}{n_A -1} \sum_{i \in A} \left( a_i - \bar{a}_A  \right)^2   + \frac{n}{n_B} \frac{1}{n_B -1}  \sum_{i \in B} \left( b_i - \bar{b}_B  \right)^2, $$}
which, under second moment conditions of potential outcomes $a$ and $b$, converges in probability to
$$ \frac{1}{p_A} \mathop {\lim}\limits_{n \rightarrow \infty } \frac{1}{n} \sum_{i=1}^n(a_i - \bar{a})^2 + \frac{1}{1-p_A} \mathop {\lim}\limits_{n \rightarrow \infty } \frac{1}{n} \sum_{i=1}^n(b_i - \bar{b})^2 . $$
Therefore, for the $\bbeta^{(a)}$ and $\bbeta^{(b)}$ defined in [\ref{beta-a-b}], the limit of $\hat \sigma^2_\textnormal{Lasso}$ is no greater than that of $\hat \sigma^2_\textnormal{unadj}$ and the difference is
{\small $$ - \mathop {\lim}\limits_{n \rightarrow \infty } \frac{1}{n}  \sum_{i=1}^{n}  \frac{1}{p_A} \biggl[ (\bx_i - \bar{\bx} )^T (\bbeta^{(a)}) \biggr]^2  +  \frac{1}{1-p_A} \biggl[ (\bx_i - \bar{\bx} )^T (\bbeta^{(b)}) \biggr]^2  .$$}
}

\noindent
{\bf Remark 7.} {\em  With the conservative variance estimate in Theorem~\ref{conservative_variance}, the Lasso adjusted confidence interval is also valid for the PATE (Population Average Treatment Effect) if there is a super population of size $N$ with $N > n$.
}

\noindent
{\bf Remark 8.} {\em The extra Condition~\ref{add-cond} is used to obtain
the following bounds for the number of selected covariates by the Lasso:
$\max{ (\hat s^{(a)},\hat s^{(b)} )} = o_p(\min{ (n_A,n_B) })$.
Condition~\ref{add-cond} can be removed from Theorem~\ref{conservative_variance}
if we redefine $\hat \sigma^2_{e^{(a)}}$ and $\hat \sigma^2_{e^{(b)}}$ without adjusting the degrees of freedom, i.e.,
\[ (\hat \sigma^*)^2_{e^{(a)}} = \frac{1}{n_A } \sum_{i \in A} \left( a_i - \bar{a}_A - (\bx_i - \bar{\bx}_A )^T \hat \bbeta^{(a)}_\textnormal{Lasso} \right)^2,  \]
\[ (\hat \sigma^*)^2_{e^{(b)}} = \frac{1}{n_B } \sum_{i \in B} \left( b_i - \bar{b}_B - (\bx_i - \bar{\bx}_B )^T \hat \bbeta^{(b)}_\textnormal{Lasso}  \right)^2, \]
and define $(\hat \sigma^*)^2_\textnormal{Lasso} = \frac{n}{n_A}  (\hat \sigma^*)^2_{e^{(a)}} + \frac{n}{n_B} (\hat \sigma^*)^2_{e^{(b)}}$. It follows from the bounds for $\max{ (\hat s^{(a)},\hat s^{(b)} )}$ that $( \hat \sigma^2_{e^{(a)}}, \hat \sigma^2_{e^{(b)}} )$ and $( (\hat \sigma^*)^2_{e^{(a)}}, (\hat \sigma^*)^2_{e^{(b)}} )$ have the same asymptotic property.}

\begin{thm}
\label{conservative_variance-2}
Assume the conditions in Theorem~\ref{normality-thm} hold. Then $(\hat \sigma^*)^2_\textnormal{Lasso}$ converges in probability to
\begin{equation*}
\frac{1}{p_A} \mathop {\lim}\limits_{n \rightarrow \infty } \sigma^2_{e^{(a)}} + \frac{1}{1-p_A} \mathop {\lim}\limits_{n \rightarrow \infty } \sigma^2_{e^{(b)}}.
\end{equation*}
\end{thm}

\noindent
{\bf Remark 9.} {\em Though $(\hat \sigma^*)^2_\textnormal{Lasso}$ has the same limit as $\hat \sigma^2_\textnormal{Lasso}$, our simulation experience shows that, in finite samples, the confidence intervals based on $(\hat \sigma^*)^2_\textnormal{Lasso}$ may yield low coverage probabilities (e.g., the coverage probability for $95\%$ confidence interval can be only $80\%$). Hence, we recommend readers to use $\hat \sigma^2_\textnormal{Lasso}$ in practice.
}


\section{Related work}
The Lasso has already made several appearances in the literature on treatment effect estimation. In the context of observational studies, \cite{Zhang2014} constructs confidence intervals for preconceived effects
or their contrasts by de-biasing the Lasso adjusted regression,
\cite{Belloni2013} employs the Lasso as a formal method for selecting adjustment variables via a two-stage procedure which concatenates features from models for treatment and outcome, and similarly, \cite{belloni2013program} gives very general results for estimating a wide range of treatment effect parameters, including the case of instrumental variables estimation. In addition to the Lasso, \cite{Li2011} considers nonparametric adjustments in the estimation of ATE. In works such as these, which deal with observational studies, confounding is the major issue. With confounding, the naive difference-in-means estimator is biased for the true treatment effect, and adjustment is used to form an unbiased estimator. However, in our work, which focuses on a randomized trial, the difference-in-means estimator is already unbiased; adjustment reduces the variance while, in fact, introducing a small amount of finite-sample bias. Another major difference between this prior work and ours is the sampling framework: we operate within the Neyman-Rubin model with fixed potential outcomes for a finite population, where the treatment group is sampled without replacement, while these papers assume independent sampling from a probability distribution with random error terms.

Our work is related to the estimation of heterogeneous or subgroup-specific treatment effects; including interaction terms to allow the imputed individual-level treatment effects to vary according to some linear combination of covariates. This is pursued in the high-dimensional setting in \cite{tian2012simple}; this work advocates solving the Lasso on a reduced set of modified covariates, rather than the full set of covariate $\times$ treatment interactions, and includes extensions to binary outcomes and survival data. The recent work in \cite{rosenblum2014optimal} considers the problem of designing multiple-testing procedures for detecting subgroup-specific treatment effects; they pose this as an optimization over testing procedures where constraints are added to enforce guarantees on type-I error rate and power to detect effects. Again, the sampling framework in these works is distinct from ours; they do not use the Neyman-Rubin model as a basis for designing the methods or investigating their properties.


\section{PAC data illustration and simulations}
We now return to the PAC-man study introduced earlier. We examine the data in more detail and explore the results of several adjustment procedures. There were 1013 patients in the PAC-man study: 506 treated (managed with PAC) and 507 control (managed without PAC, but retaining the option of using   alternative devices). The outcome variable is quality-adjusted life years (QALYs). One QALY represents one year of life in full health; in-hospital death corresponds to a QALY of zero. We have $59$ covariates about each individual in the study; we include all main effects as well as $1113$ two-way interactions, and form a design matrix $\mathbf{X}$ with $1172$ columns and $1013$ rows. See Appendix B for more details on the design matrix.


\ABedit{The assumptions that underpin the theoretical guarantees of the $\ATElasso$ estimator are, in practice, not explicitly checkable, but we attempt to inspect the quantities that are involved in the conditions to help readers make their own judgement. The uniform bounds on the fourth moments refer to a hypothetical sequence of populations; these cannot be verified given that the investigator has a single dataset.  However, as an approximation, the fourth moments of the data can be inspected to ensure that they are not too large. In this data set, the maximum fourth moment of the covariates is $37.3$, \ABedit{which is indicative of a heavy-tailed and potentially destabilizing covariate; however, it} occurs in an interaction term not selected by the lasso, and thus does not influence the estimate\footnote{The fourth moments of the covariates are shown in Fugure~\ref{fig:fourth-moment} in Appendix F. The covariates with the largest two fourth moments ($37.3$ and $34.9$ respectively) are quadratic term $interactnew^2$ and interaction term $IMscorerct:systemnew$. Neither of them are selected by the Lasso to do the adjustment.}. Checking the conditions for high-dimensional consistency of the Lasso would require knowledge of the unknown active set $S$, and moreover, even if it were known, calculating the cone invertibility factor would involve an infeasible optimization.  This is a general issue in the theory of sparse linear high-dimensional estimation.  To approximate these conditions, we use the bootstrap to estimate the active set of covariates $S$ and the error terms $e^{(a)}$ and $e^{(b)}$. See the supporting information for more details. Our estimated $S$ contains $16$ covariates and the estimated second moments of $e^{(a)}$ and $e^{(b)}$ are $11.8$ and $12.0$, respectively. The estimated maximal covariance $\delta_n$ equals $0.34$ and the scaling $(s \log p)/{\sqrt n} $ is $3.55$. \ABedit{While this is not close to zero,} we should mention that the estimation of $\delta_n$ and $(s \log p)/{\sqrt n}$ can be unstable and less accurate since it is based on a subsample of the population. As an approximation to Condition \ref{last-cond}, we examine the largest and smallest eigenvalues of the sub-Gram matrix $(1/n)\mathbf{X}_S^T\mathbf{X}_S$, which are $2.09$ and $0.18$ respectively. Thus the quantity in Condition \ref{last-cond} seems reasonably bounded away from zero.}

We now estimate the ATE using the unadjusted estimator, the Lasso adjusted estimator and the OLS adjusted estimator which is computed based on a sub-design matrix containing only the 59 main effects. We also present results for the two-step estimator $\widehat{ATE}_\textnormal{Lasso+OLS}$ which adopts the Lasso to select covariates and then uses OLS to refit the regression coefficients. See \cite{Efron2004,Meinshausen2007,Belloni2009LassoOLS,Liu2013} for statistical properties of Lasso+OLS estimator in linear regression model. Let $\hat \bbeta^{(a)}$ be the Lasso estimator defined in \ref{def-lasso-a} (we omit the subscript ``Lasso" for the sake of simplicity) and let $\hat S^{(a)} = \{ j: \hat \bbeta^{(a)}_j \neq 0 \}$ be the support of $\hat \bbeta^{(a)}$. The Lasso+OLS adjustment vector $\hat \beta^{(a)}_\textnormal{Lasso+OLS}$ for treatment group A is defined by
\begin{equation}
\begin{split}
\label{def-lassools-a}
\hat \bbeta^{(a)}_\textnormal{Lasso+OLS}  = \argmin_{\bbeta: \ \beta_j=0, \ \forall j \notin \hat S^{(a)} } & \frac{1}{2n_A}\sum_{i\in A} \left[ a_i - \bar{a}_A \right. \\
&   \left. - ( \bx_i - \bar{\bx}_A )^T \bbeta \right]^2. \nonumber
\end{split}
\end{equation}
We can define the Lasso+OLS adjustment vector $\hat \bbeta^{(b)}_\textnormal{Lasso+OLS}$ for control group B similarly. Then $\widehat{ATE}_\textnormal{Lasso+OLS}$ is given by
\begin{equation}
\begin{split}
\widehat{ATE}_\textnormal{Lasso+OLS} = & \left[ \bar{a}_A - ( \bar{\bx}_A - \bar{\bx} )^T \hat \bbeta^{(a)}_\textnormal{Lasso+OLS} \right] \\
& - \left[ \bar{b}_B - ( \bar{\bx}_B - \bar{\bx} )^T \hat \bbeta^{(b)}_\textnormal{Lasso+OLS} \right]. \nonumber
\end{split}
\end{equation}
In the next paragraph and in Algorithm~\ref{alg:cv} of Appendix F, we show how we adapt the cross-validation procedure to select the tuning parameter for $\widehat{ATE}_\textnormal{Lasso+OLS}$ based on a combined performance of Lasso and OLS, or cv(Lasso+OLS).

We use the R package ``glmnet" to compute the Lasso solution path and select the tuning parameters $\lambda_a$ and $\lambda_b$ by 10-fold Cross Validation (CV). To indicate the method of selecting tuning parameters, we denote the corresponding estimators as cv(Lasso) and cv(Lasso+OLS) respectively. We should mention that for the cv(Lasso+OLS) adjusted estimator, we compute the CV error for a given value of $\lambda_a$ (or $\lambda_b$) based on the whole Lasso+OLS procedure instead of just the Lasso estimator (see Algorithm~\ref{alg:cv} in Appendix F). Therefore, the cv(Lasso+OLS) and the cv(Lasso) may select different covariates to do the adjustment. This type of cross validation requires more computation than the cross validation based on just the Lasso estimator since it needs to compute the OLS estimator for each fold and each given $\lambda_a$ (or $\lambda_b$), but it can give better prediction and model selection performance.

\begin{figure}
\captionsetup{width=.45\textwidth}
\centerline{\includegraphics[width=.45\textwidth, height=.45\textwidth]{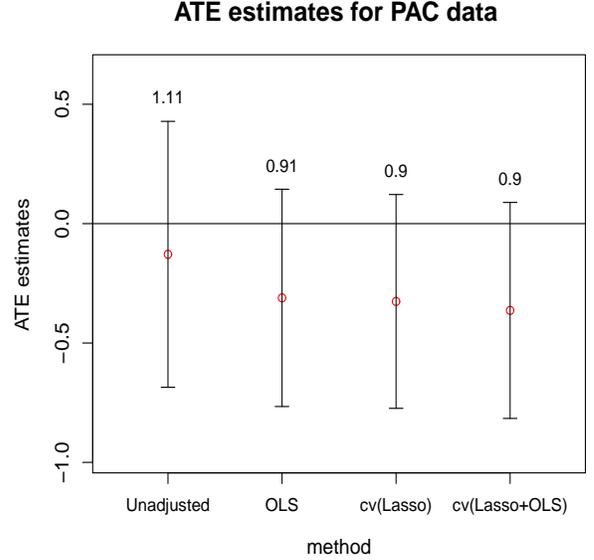}}
\caption{ATE estimates (red circles) and $95\%$ confidence intervals (bars) for the PAC data. The numbers above each bar are the corresponding interval lengths.}\label{fig:pac}
\end{figure}

Figure~\ref{fig:pac} presents the ATE estimates along with $95\%$
confidence intervals (CI). The interval lengths are shown on top of
each interval bar. All the methods give confidence intervals
containing $0$; hence, this experiment failed to provide sufficient evidence to reject the hypothesis that PAC did not have an effect on patient QALYs (either positive or negative). Since the caretakers of patients managed without PAC retained the option of using less invasive cardiac output monitoring devices, such an effect may have been particularly hard to detect in this experiment.

\begin{table*}[!ht]
\begin{center}
 \caption{\label{tab:pac} Statistics for the PAC illustration}
\begin{tabular*}{\hsize}{@{\extracolsep{\fill}}lccccc}
  & & & & & No. of selected covariates  \\
  \cline{5-6}
 Methods & $\widehat{ATE}$ & $\hat \sigma_\textnormal{ATE}$ & $95\% $ confidence interval & treated & control \\
  \hline
Unadjusted & -0.13 & 0.081 & [-0.69,0.43] & - & - \\
  OLS & -0.31 & 0.054 & [-0.77,0.14] & - & - \\
  cv(Lasso) & -0.33 & 0.052 & [-0.77,0.12] & 24 & 8 \\
  cv(Lasso+OLS) & -0.36 & 0.053 & [-0.82,0.09] & 4 & 5 \\
   \hline
\end{tabular*}
\end{center}
\end{table*}

However, it is interesting to note that (see Table~\ref{tab:pac}), compared with the unadjusted estimator, the OLS adjusted estimator causes the ATE estimate to decrease (from -0.13 to -0.31), and shortens the confidence interval by about $20\%$. This is due mainly to the imbalance in the pre-treatment probability of death, which was highly predictive of the post-treatment QALYs. The cv(Lasso) adjusted estimator yields a comparable ATE estimate and confidence interval, but the fitted model is more interpretable and parsimonious than the OLS model: it selects $24$ and $8$ covariates for treated and control, respectively. The cv(Lasso+OLS) estimator selects even fewer covariates: $4$ and $5$ for treated and control, respectively, but performs a similar adjustment as the cv(Lasso) (see the comparison of fitted values in Figure~\ref{fig:fitvalue}). We also note that these adjustments agree with the one performed in \cite{Miratrix2013}, where the treatment effect was adjusted downwards to $-0.27$ after stratifying into 4 groups based on predicted probability of death.

The covariates selected by Lasso for adjustment are shown in Table~\ref{tab:selected-covariate}, where ``A$\cdot$A" denote quadratic term of the covariate A and ``A:B" denote two way interaction between two covariates A and B. Among them, patient's age and estimated probability of death (p\_death), together with the quadratic term ``age$\cdot$age" and interactions ``age:p\_death" and ``p\_death:mech\_vent\footnote{mechanical ventilation at admission}", are the most important covariates for the adjustment. The patients in control group are slightly older and have slightly higher risk of death. These covariates are important predictors of the outcome. Therefore, the unadjusted estimator \ABedit{may overestimate} the benefits of receiving PAC.

\begin{table*}[ht]
 \caption{\label{tab:selected-covariate} Selected covariates for adjustment}
\begin{tabular*}{\hsize}{@{\extracolsep{\fill}\vrule height 10.5pt depth4pt width0pt}lcc}
method & treatment & covariates \\ \hline
cv(Lasso+OLS) & treated & age, p\_death, age$\cdot$age, age:p\_death \\ \hline

cv(Lasso+OLS) & control & age, p\_death, age$\cdot$age, age:p\_death, p\_death:mech\_vent \\ \hline

cv(Lasso) & treated & pac\_rate, age, p\_death, age$\cdot$age, p\_death$\cdot$p\_death, region:im\_score, region:systemnew, \\
 & &  pac\_rate:age, pac\_rate:p\_death, pac\_rate:systemnew, im\_score:interactnew, age:p\_death,  \\
& &  age:glasgow, age:systemnew, interactnew:systemnew, pac\_rate:creatinine,  \\
& & age:mech\_vent, age:respiratory, age:creatinine, interactnew:mech\_vent,  \\
& & interactnew:male, glasgow:organ\_failure, p\_death:mech\_vent, systemnew:male \\ \hline

cv(Lasso) & control & age, p\_death, age$\cdot$age, unitsize:p\_death, pac\_rate:systemnew, age:p\_death, \\
& &  interactnew:mech\_vent, p\_death:mech\_vent \\ \hline
\end{tabular*}
Covariate definitions: age (patient's age); p\_death (baseline probability of death); mech\_vent (mechanical ventilation at admission); region (geographic region); pac\_rate (PAC rate in unit); creatinine, respiratory, glasgow, interactnew, organ\_failure, systemnew, im\_score (various physiological indicators).
\end{table*}

Since not all the potential outcomes are observed, we cannot know the true gains of adjustment methods. However, we can estimate the gains via building a simulated set of potential outcomes by matching treated units to control units on observed covariates. We use the matching method described in \cite{Sekhon2005} which gives $1013$ observations with all potential outcomes \ABedit{imputed}. 
We match on the $59$ main effects \ABedit{only}. The ATE is $-0.29$. We then use this synthetic data set to calculate the biases, standard deviations (SD) and root-mean square errors ($\sqrt{\textnormal{MSE}}$) of different ATE estimators based on $25000$ replicates of completely randomized experiment which assigns 506 subjects to the treated group and the remainders to the control group.

\begin{figure}
\captionsetup{width=.45\textwidth}
\centerline{\includegraphics[width=0.45\textwidth]{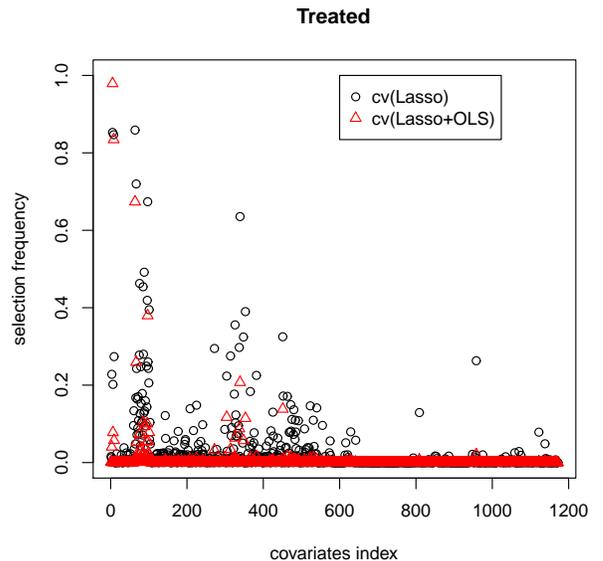}}
\caption{Selection stability comparison of cv(Lasso) and cv(Lasso+OLS) for treatment group.}\label{fig:stability}
\end{figure}

\begin{figure}
\captionsetup{width=.45\textwidth}
\centerline{\includegraphics[width=0.45\textwidth]{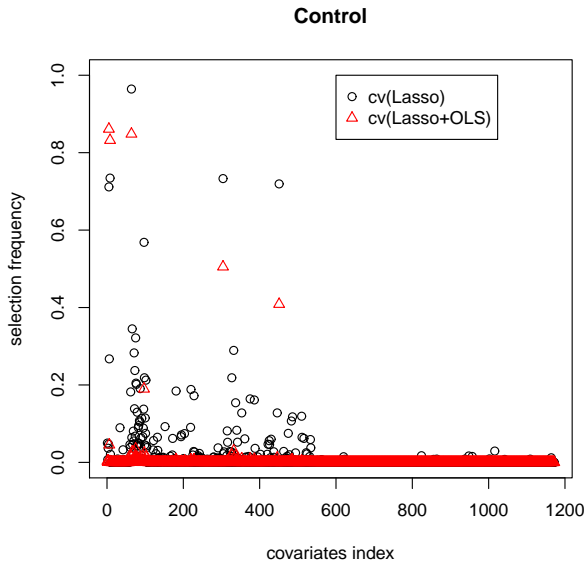}}
\caption{Selection stability comparison of cv(Lasso) and cv(Lasso+OLS) for control group.}\label{fig:stability2}
\end{figure}

Table~\ref{tab:pac-simu} shows the results. For all the methods, the bias is substantially smaller (\ABedit{by a factor of 100}) than the SD. The SD and $\sqrt{ \textnormal{MSE} }$ of the OLS adjusted estimator are both $10.2\%$ smaller than those of the unadjusted estimator, while the cv(Lasso) and cv(Lasso+OLS) adjusted estimators further improve the SD and $\sqrt{ \textnormal{MSE} }$ of the OLS adjusted estimator by approximately $4.7\%$. Moreover, all these methods provide conservative confidence intervals with coverage probabilities higher than $99\%$. However, the interval lengths of the OLS, cv(Lasso) and cv(Lasso+OLS) adjusted estimator are comparable and are approximately $10\%$ shorter than that of the unadjusted estimator. The cv(Lasso+OLS) adjusted estimator is similar to the cv(Lasso) adjusted estimator in terms of mean squared error, confidence interval length and coverage probability, but outperforms the latter with much fewer and more stable covariates in the adjustment (see Figures~\ref{fig:stability} and ~\ref{fig:stability2} for the selection frequency of each covariate for treatment group and control group respectively). We also show in Figure~\ref{fig:density} that the sampling distribution of the estimates is very close to Normal.

We conduct additional simulation studies to evaluate the finite sample performance of $\widehat{ATE}_\textnormal{Lasso}$ and compare it with that of the OLS adjusted estimator and the unadjusted estimator. A qualitative analysis of these simulations yields the same conclusions as presented above; however, for the sake of brevity, we defer the simulation details in the supporting information.

\begin{table*}
 \caption{\label{tab:pac-simu} Statistics for the PAC synthetic data set}
\begin{tabular*}{\hsize}{@{\extracolsep{\fill}}lccccccc}
    & & & & & & \multicolumn{2}{c}{No. of selected covariates}  \\
  \cline{7-8}
 & Bias & SD & $\sqrt{ \textnormal{MSE} }$ & Coverage $(\%)$ & Length & treated  & control \\ \hline

  Unadjusted       & 0.001(0) & 0.20(0.02) & 0.20(0.02) & 99 & 1.06 & - & - \\
  OLS              & 0.002(0) & 0.18(0.02) & 0.18(0.02) & 99 & 0.95 & - & - \\
  cv(Lasso)        & 0.001(0) & {\bf 0.17(0.02)} & {\bf 0.17(0.02)} & 99 & {\bf 0.94} & 25(23) & 15(14) \\
  cv(Lasso+OLS)    & {\bf 0.000(0)} & {\bf 0.17(0.02)} & {\bf 0.17(0.02)} & 99 & 0.95 & {\bf 6(6)} & {\bf 4(3)} \\  \hline
\end{tabular*}
The numbers in parentheses are the corresponding standard errors estimated by using the bootstrap with $B=500$ resamplings of the ATE estimates.
\end{table*}



\section{Discussion}
We study the Lasso adjusted average treatment effect (ATE) estimate under the Neyman-Rubin model for randomization. Our purpose in using the Neyman-Rubin model was to investigate the performance of the Lasso under a realistic sampling framework which does not \JSedit{impose} strong assumptions on the data. We provide conditions that ensure asymptotic normality, and provide a Neyman-type estimate of the asymptotic variance which can be used to construct a conservative confidence interval for the ATE. While we do not require an explicit generative linear model to hold, our theoretical analysis requires the existence of latent `adjustment vectors' such that moment conditions of the error terms are satisfied, and \JSedit{that} the cone invertibility condition of the sample covariance matrix is satisfied in addition
to moment conditions for OLS adjustment as in \cite{Lin2013}. Both assumptions are difficult to check in practice.  \ABedit{In our theory, we do not address whether these assumptions are necessary for our results to hold, though simulations indicate that the moment conditions cannot be substantially weakened.}  As a by-product of our analysis, we extend Massart's concentration inequality for sampling without replacement, which is useful for theoretical analysis under the Neyman-Rubin model. Simulation studies and \JSedit{the} real data illustration show the advantage of the Lasso-adjusted estimator in terms of estimation accuracy and model interpretation. In practice, we recommend a variant of Lasso, cv(Lasso+OLS), to select covariates and perform the adjustment, since it gives similar coverage probability and confidence interval length when compared with cv(Lasso), but with far fewer covariates selected. In future work, we plan to extend our analysis to other popular methods in high-dimensional statistics such as Elastic-Net and ridge regression, which may be more appropriate for estimating adjusted ATE under different assumptions.

The main goal of using Lasso in this paper is to reduce the variance (and overall mean squared error) of ATE estimation. \JSedit{Another important task is to estimate heterogenous treatment effects and provide conditional treatment effect estimates for subpopulations. When the Lasso models of treatment and control outcomes are different, both in variables selected} and coefficient values, this could be interpreted as modeling treatment effect heterogeneity in terms of covariates. \JSedit{However, reducing variance of the ATE estimate and estimating heterogenous treatment effects have completely different targets.  Targeting heterogenous treatment effects may result in more variable ATE estimates.}  Moreover, our simulations show that the set of covariates selected by the Lasso is unstable and this may cause problems when interpreting them as evidence of heterogenous treatment effects. \JSedit{How best to estimate such effects is an open question that we would like to study in future research.}

\section{Materials and Methods}
We did not conduct the PAC-man experiment, and we are analyzing secondary data without any personal identifying information. As such, this study is exempt from human subjects review. The original experiments underwent human subjects review in the UK~\cite{harvey2005assessment}.

\section*{Acknowledgements}
We thank David Goldberg for helpful discussions, Rebecca Barter for copyediting and suggestions for clarifying the text, and Winston Lin for comments. We thank Richard Grieve (LSHTM), Sheila Harvey (LSHTM), David Harrison (ICNARC) and Kathy Rowan (ICNARC) for access to data from the PAC-Man CEA and the ICNARC CMP database. This research is partially supported by NSF grants DMS-11-06753, DMS-12-09014, DMS-1107000, DMS-1129626, DMS-1209014, CDS$\&$E-MSS, 1228246DMS-1160319 (FRG), AFOSR grant FA9550-14-1-0016, NSA Grant H98230-15-1-0040, the Center for Science of Information (CSoI), an US NSF Science and Technology Center, under grant agreement CCF-0939370, the Department of Defense (DoD) for Office of Naval Research (ONR) grant N00014-15-1-2367 and the National Defense Science \& Engineering Graduate Fellowship (NDSEG) Program.

\appendix

\section{Simulation}

In this section we carry out simulation studies to evaluate the finite sample performance of $\widehat{ATE}_\textnormal{Lasso}$ estimator. We also present results for the $\widehat{ATE}_\textnormal{OLS}$ estimator when $p<n$ and the two-step estimator $\widehat{ATE}_\textnormal{Lasso+OLS}$.

We use the R package ``glmnet" to compute the Lasso solution path. We select the tuning parameters $\lambda_a$ and $\lambda_b$ by 10-fold Cross Validation (CV) and denote the corresponding adjusted estimators as cv(Lasso) and cv(Lasso+OLS) respectively. We should mention that for the cv(Lasso+OLS) adjusted estimator, we compute the CV error for a given value of the $\lambda_a$ (or $\lambda_b$) based on the whole Lasso+OLS estimator instead of the Lasso estimator, see Algorithm~\ref{alg:cv} for details. Therefore, the cv(Lasso+OLS) adjusted estimator and the cv(Lasso) adjusted estimator may select different covariates to do the adjustment. This type of cross validation for cv(Lasso+OLS) requires more computation effort than the cross validation based on just the Lasso estimator since it needs to compute the OLS estimator for each fold and for each $\lambda_a$ (or $\lambda_b$), but it can give better prediction and covariates selection performance.

The potential outcomes $a_i$ and $b_i$ are generated from the following nonlinear model: for $i=1,...,n$,
\[ a_i = \sum_{j=1}^s x_{ij} \beta_{j}^{(a1)} + \exp{ \left( \sum_{j=1}^s x_{ij} \beta_{j}^{(a2)} \right) } + \epsilon^{(a)}_i, \]
\[ b_i = \sum_{j=1}^s x_{ij} \beta_{j}^{(b1)} + \exp{ \left( \sum_{j=1}^s x_{ij} \beta_{j}^{(b2)} \right) } + \epsilon^{(b)}_i, \]
where $\epsilon^{(a)}_i$ and $\epsilon^{(b)}_i$ are independent error terms. We set $n=250$, $s=10$, $p=50$ and $500$. For $p=50$, we can compute OLS estimator and compare it with the Lasso. The covariates vector $\bx_i$ is generated from a multivariate normal distribution $\mathcal{N}(0,\Sigma)$. We consider two different Toeplitz covariance matrices $\Sigma$ which control the correlation among the covariates:
\[ \Sigma_{ii}=1; \   \Sigma_{ij}=\rho^{|i-j|} \ \forall i\neq j, \]
where $\rho=0,0.6$. The true coefficients $\beta^{(a1)}_j$, $\beta^{(a2)}_j$, $\beta^{(b1)}_j$, $\beta^{(b2)}_j $ are generated independently according to
\[ \beta_{j}^{(a1)} \sim  t_3; \ \ \   \beta_{j}^{(a2)} \sim 0.1*t_3,  \ \ \ j=1,...,s, \]
\[ \beta_{j}^{(b1)} \sim \beta_{j}^{(a1)} + t_3; \ \ \   \beta_{j}^{(b2)} \sim \beta_{j}^{(a2)} + 0.1*t_3,  \ \ \ j=1,...,s, \]
where $t_3$ denotes the $t$ distribution with three degrees of freedom. This ensures that the treatment effects are not not constant across individuals, and that the linear model does not hold in this simulation. The error terms $\epsilon^{(a)}_i$ and $\epsilon^{(b)}_i$ are generated according to the following linear model with some hidden covariates $\bz_i$:
\[ \epsilon^{(a)}_i = \sum_{j=1}^s z_{ij} \beta_{j}^{(a1)} +  \tilde \epsilon^{(a)}_i, \]
\[ \epsilon^{(b)}_i = \sum_{j=1}^s z_{ij} \beta_{j}^{(b1)} +  \tilde \epsilon^{(b)}_i, \]
where $\tilde \epsilon^{(a)}_i$ and $\tilde \epsilon^{(b)}_i$ are drawn independently from standard normal distribution. The vector $\bz_i$ is independent of $\bx_i$ and also drawn independently from the multivariate normal distribution $\mathcal{N}(0,\Sigma)$. The values of $\bx_i$, $\beta^{(a1)}$, $\beta^{(a2)}$, $\beta^{(b1)}$, $\beta^{(b2)} $, $\bz_i$, $\tilde \epsilon^{(a)}_i$, $\tilde \epsilon^{(b)}_i$, $a_i$ and $b_i$ are generated once and then kept fixed.

After the potential outcomes are generated, a completely randomized experiment is simulated $25000$ times, assigning $n_A= 100,125,150$ subjects to treatment A and the remainder to control B. There are $12$ different combinations of $(p,\rho,n_A)$ in total.

\begin{figure*}
\centerline{\includegraphics[width=.75\textwidth]{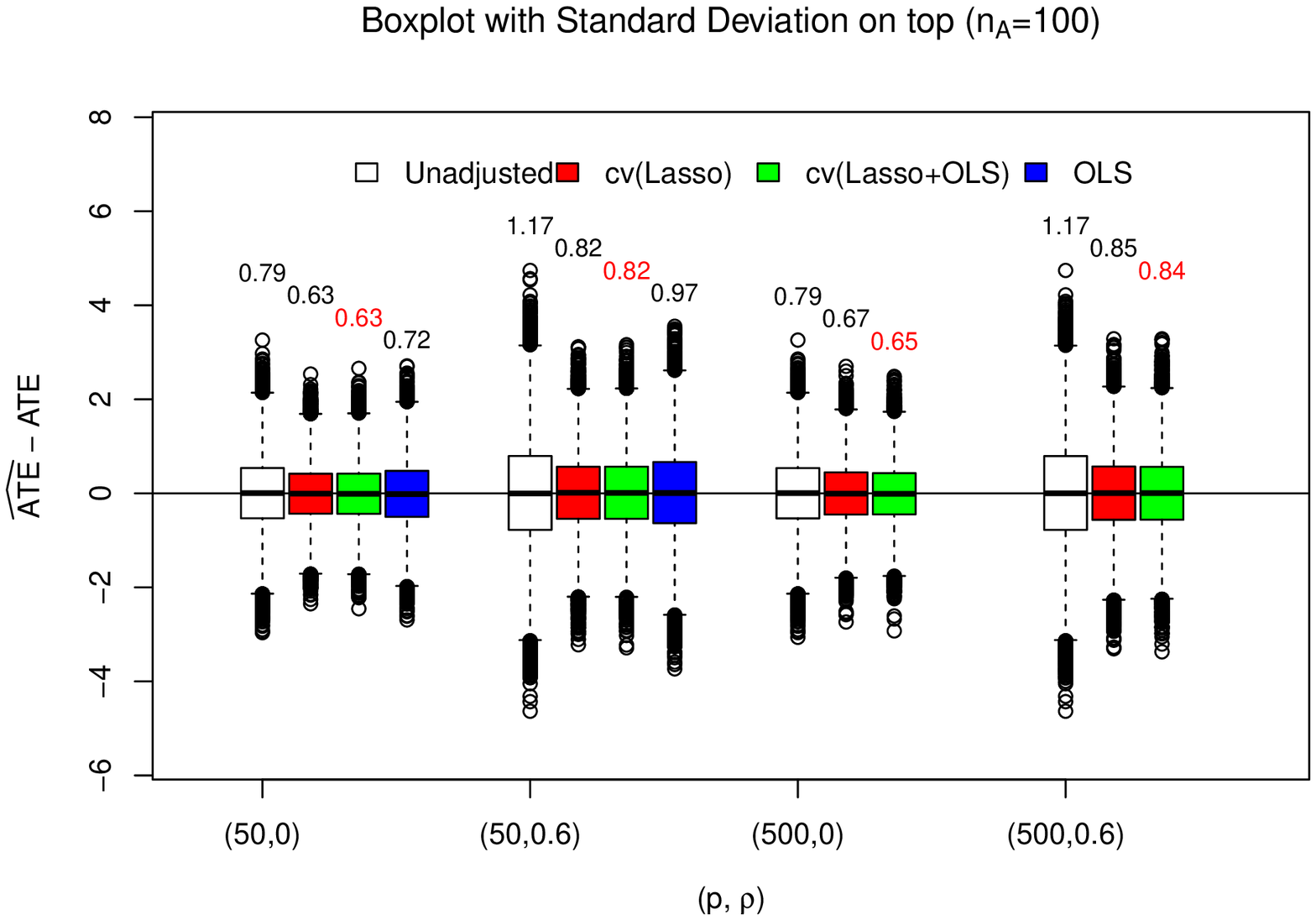}}
\vspace*{0.05in}
\caption{Boxplot of the unadjusted, OLS adjusted (only computed when $p=50$), cv(Lasso) and cv(Lasso+OLS) adjusted estimators with their standard deviations presented on top of each box for $n_A=100$.}\label{fig:boxplot1}
\end{figure*}
\begin{figure*}
\centerline{\includegraphics[width=.75\textwidth]{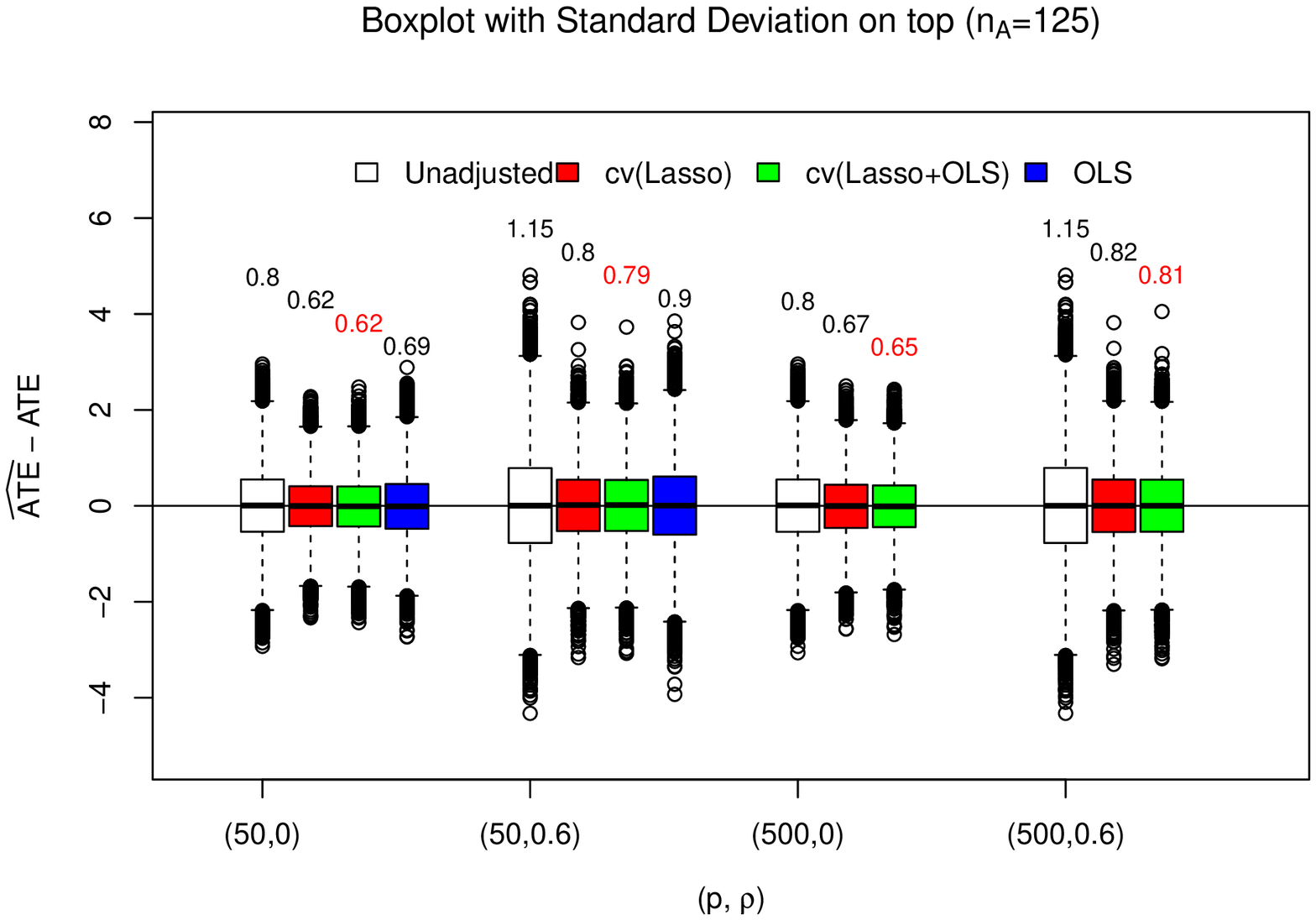}}
\vspace*{0.05in}
\caption{Boxplot of the unadjusted, OLS adjusted (only computed when $p=50$), cv(Lasso) and cv(Lasso+OLS) adjusted estimators with their standard deviations presented on top of each box for $n_A=125$.}\label{fig:boxplot2}
\end{figure*}
\begin{figure*}
\centerline{\includegraphics[width=.75\textwidth]{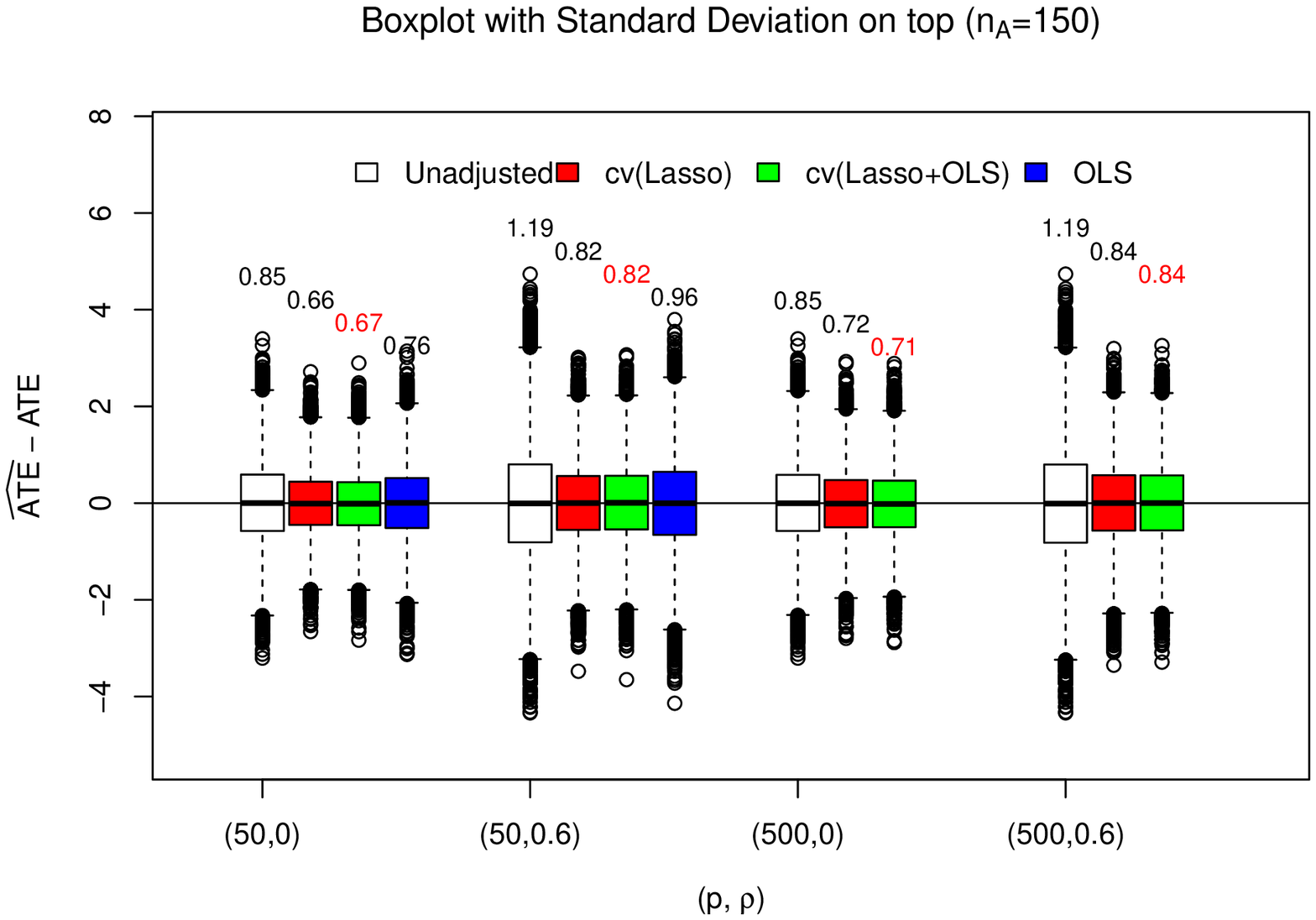}}
\vspace*{0.05in}
\caption{Boxplot of the unadjusted, OLS adjusted (only computed when $p=50$), cv(Lasso) and cv(Lasso+OLS) adjusted estimators with their standard deviations presented on top of each box for $n_A=150$.}\label{fig:boxplot3}
\end{figure*}

Figures~\ref{fig:boxplot1}, \ref{fig:boxplot2}, \ref{fig:boxplot3} show the boxplot of different ATE estimators with their standard deviations (computed from $25000$ replicates of randomized experiments) presented on top of each box. Regardless of whether the design is balanced $(n_A=125)$ or not $(n_A=100,150)$, the regression based estimators have much smaller variances and than that of the unadjusted estimator and therefore improve the estimation precision.


To further compare the performance of these estimators, we present the bias, the standard deviation (SD) and the root-mean square error ($\sqrt{\textnormal{MSE}}$) of the estimates in Table~\ref{tab:mse}. Bias is reported as the absolute difference from the true treatment effect. We find that the bias of each method is substantially smaller (more than 10 times smaller) than the SD. The cv(Lasso) and cv(Lasso+OLS) adjusted estimators perform similar in terms of SD and $\sqrt{\textnormal{MSE}}$: reducing those of the OLS adjusted estimator and the unadjusted estimator by $10\% - 15\%$ and $15\% - 31\%$ respectively. We also compare the number of selected covariates by cv(Lasso) and cv(Lasso+OLS) for treatment group and control group separately, see Table \ref{tab:modelsize}. It is easy to see that the cv(Lasso+OLS) adjusted estimator uses many fewer (more than $44\%$) covariates in the adjustment to obtain similar improvement of SD and $\sqrt{\textnormal{MSE}}$ of ATE estimate as the cv(Lasso) adjusted estimator. Moreover, we find that the covariates selected by the cv(Lasso+OLS) are more stable across different realizations of treatment assignment than the covariates selected by the cv(Lasso). Overall, the cv(Lasso+OLS) adjusted, the cv(Lasso) adjusted, the OLS adjusted and the unadjusted estimators perform from best to worst.



\begin{figure*}
\centerline{\includegraphics[width=0.7\textwidth]{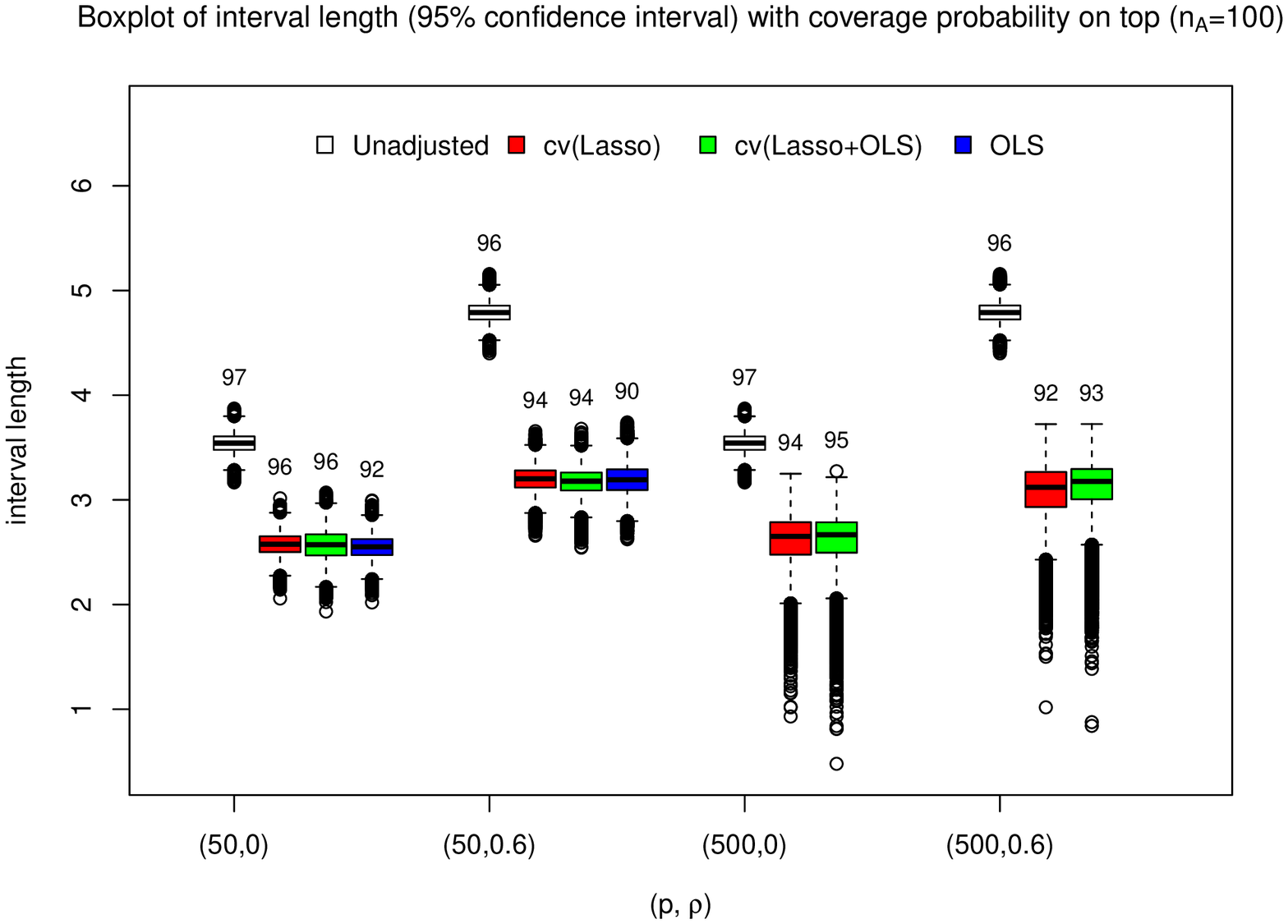}}
\vspace*{0.05in}
\caption{Boxplot of the interval length with coverage probability $(\%)$ on top of each box for the unadjusted, OLS adjusted (only computed when $p=50$), cv(Lasso) adjusted and cv(Lasso+OLS) adjusted estimators with $n_A=100$.}\label{fig:boxplot-interval1}
\end{figure*}

\begin{figure*}
\centerline{\includegraphics[width=0.7\textwidth]{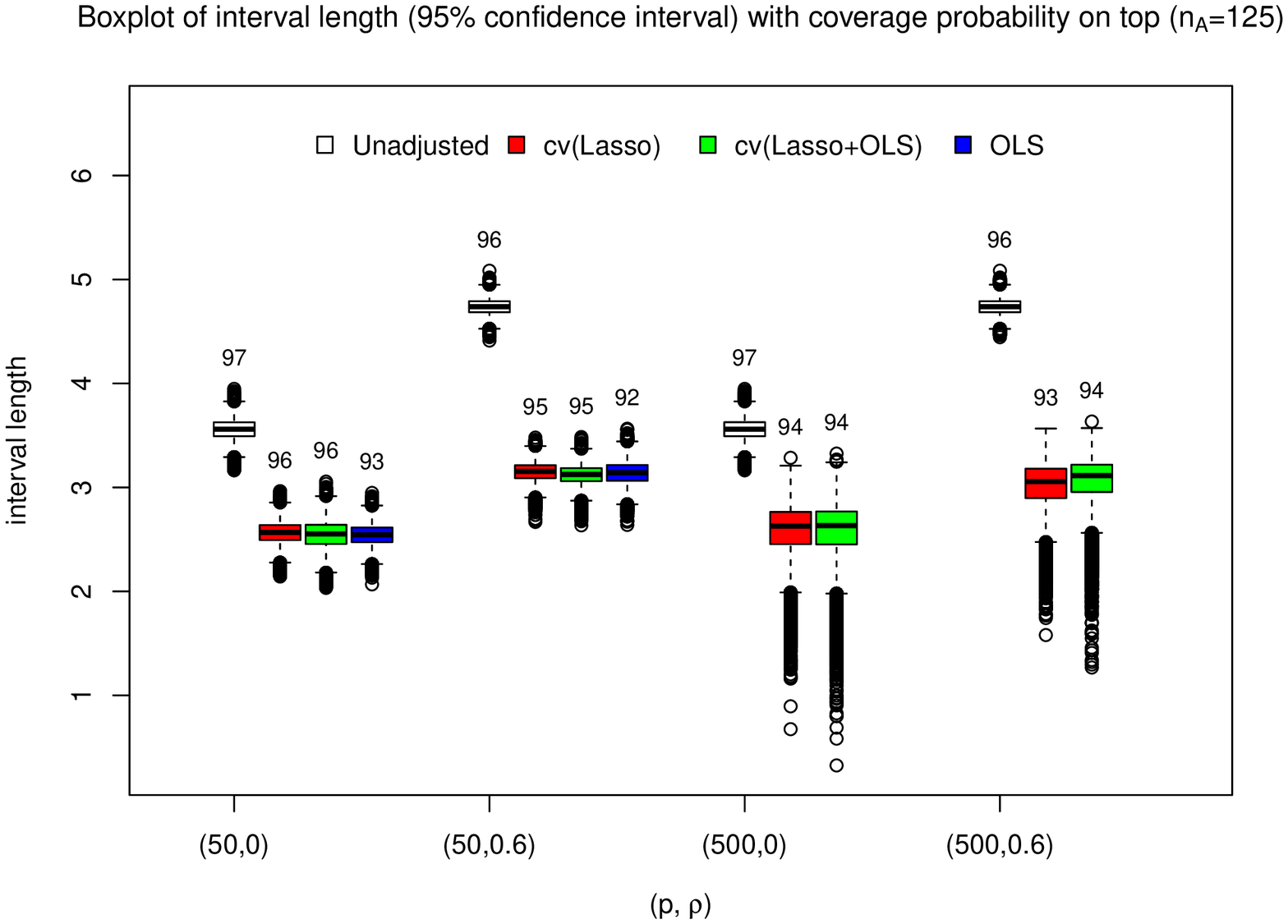}}
\vspace*{0.05in}
\caption{Boxplot of the interval length with coverage probability $(\%)$ on top of each box for the unadjusted, OLS adjusted (only computed when $p=50$), cv(Lasso) adjusted and cv(Lasso+OLS) adjusted estimators with $n_A=125$.}\label{fig:boxplot-interval2}
\end{figure*}

\begin{figure*}
\centerline{\includegraphics[width=0.7\textwidth]{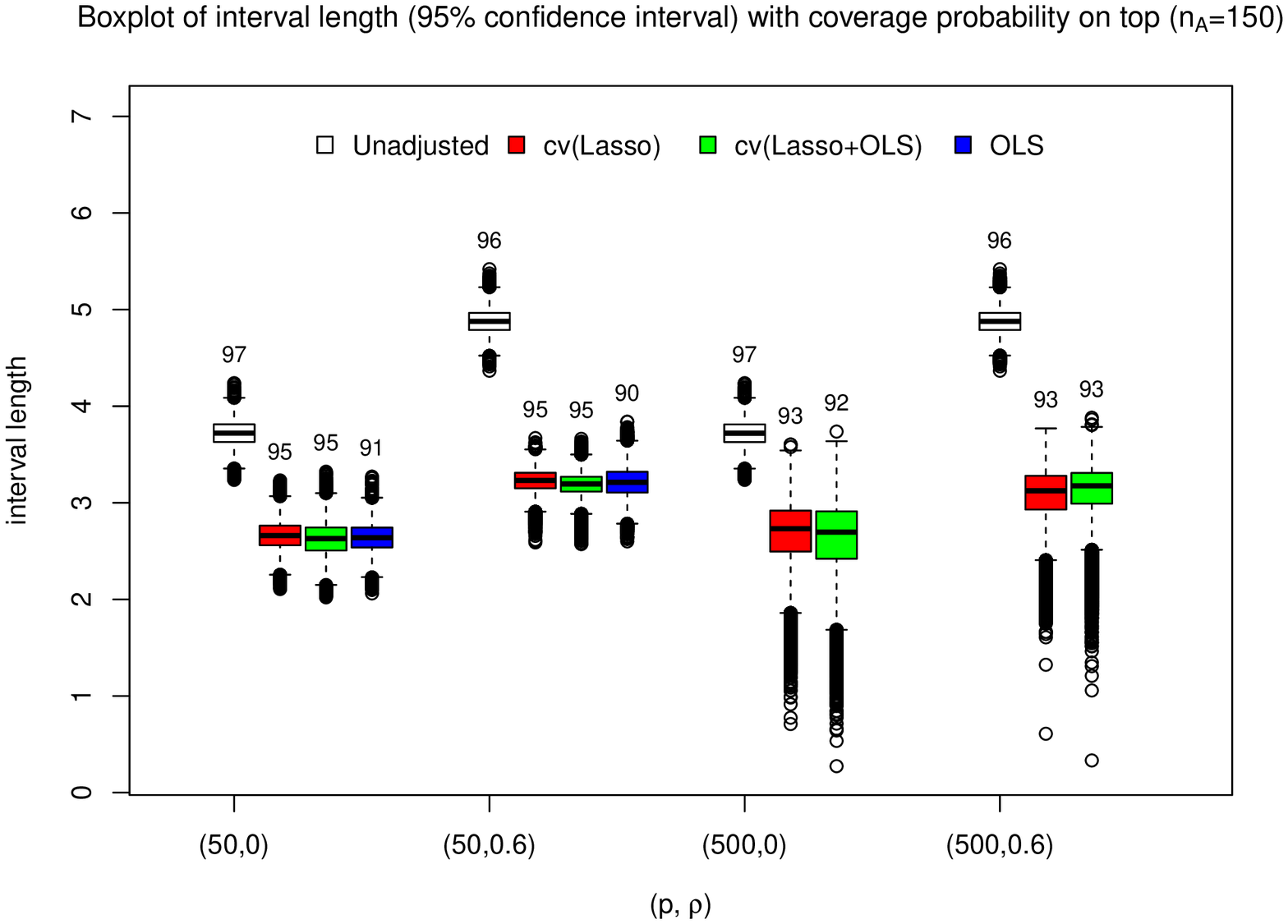}}
\vspace*{0.05in}
\caption{Boxplot of the interval length with coverage probability $(\%)$ on top of each box for the unadjusted, OLS adjusted (only computed when $p=50$), cv(Lasso) adjusted and cv(Lasso+OLS) adjusted estimators with $n_A=150$.}\label{fig:boxplot-interval3}
\end{figure*}

We move now to study the finite sample performance of Neyman-type
conservative variance estimates. For each simulation example and each
one of the $25000$ completely randomized experiments, we calculate the
ATE estimates ($\widehat{ATE}$) and the Neyman variance estimates
($\hat \sigma$) and then form the $95\%$ confidence intervals
$[\widehat{ATE} - 1.96\cdot \hat \sigma/ \sqrt{n}, \widehat{ATE} +
1.96\cdot \hat \sigma/ \sqrt{n}]$. Figures~\ref{fig:boxplot-interval1}, \ref{fig:boxplot-interval2}, \ref{fig:boxplot-interval3}  present
the boxplot of the interval length with the coverage probability noted
on top of each box for the unadjusted, OLS adjusted (only computed
when $p=50$), cv(Lasso) adjusted and cv(Lasso+OLS) adjusted estimators. More
results are showed in Table \ref{tab:coverage}. We find that all
the confidence intervals for the unadjusted estimator are conservative. The cv(Lasso) adjusted and the cv(Lasso+OLS0 adjusted estimators perform very similar:
although their coverage probability (at least $92\%$) may be slightly less than the pre-assigned confidence level ($95\%$), their mean interval length is much shorter ($26\% - 37\%$) than that of the unadjusted estimator. The OLS adjusted estimator has comparable interval length with the cv(Lasso) and cv(Lasso+OLS) adjusted estimator, but has slightly worse coverage probability ($90\% - 93\%$).

\begin{figure*}
\centering\includegraphics[width=0.9\textwidth]{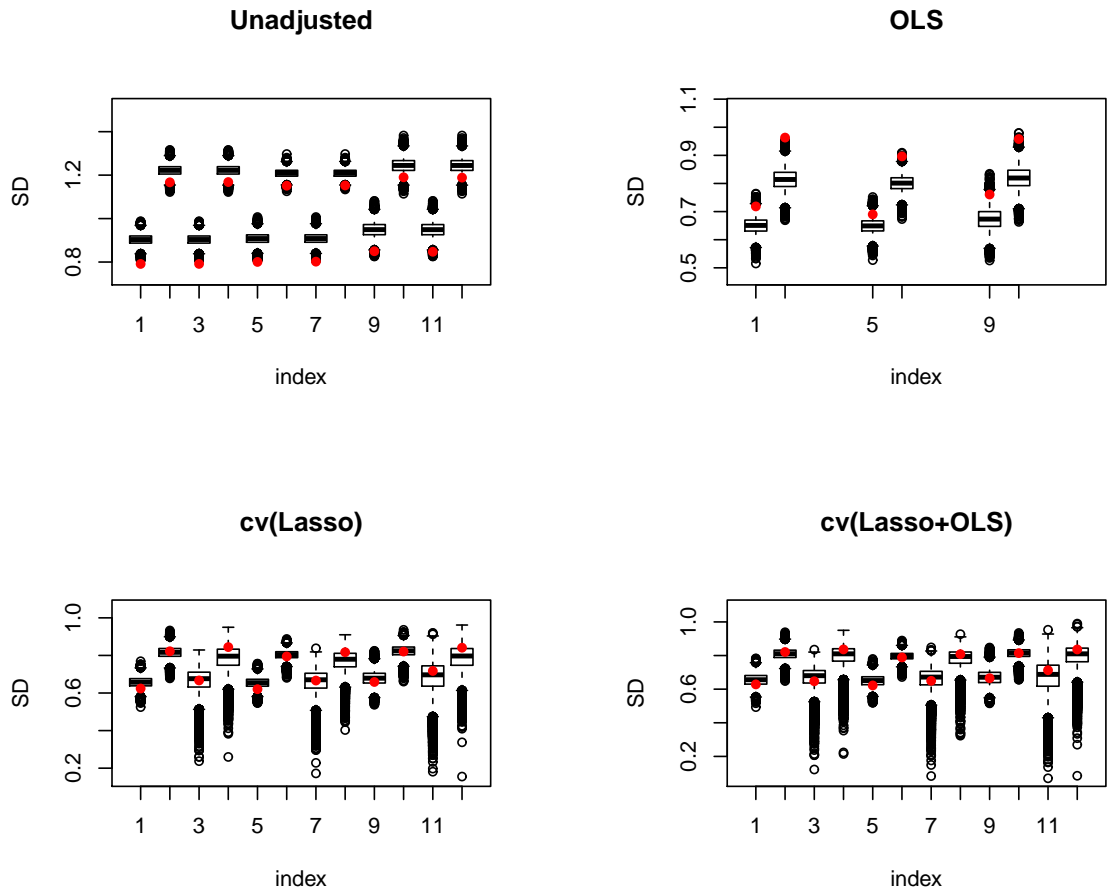}
\vspace*{0.05in}
\caption{Boxplot of Neyman SD estimate with the ``true" SD presented as red dot.}
\label{fig:com-var-neyman}
\end{figure*}

To further investigate how good the Neyman standard deviation (SD) estimate is, we compare them in Figure~\ref{fig:com-var-neyman} with the ``true'' SD presented in Table~\ref{tab:mse} (the SD of the ATE estimates over $25000$ randomized experiments). We find that Neyman SD estimate is very conservative for the unadjusted estimator (its mean is $5\% - 14\%$ larger than the ``true" SD); while for the OLS adjusted estimator, the mean of Neyman SD estimate can be $6\% - 100\%$ smaller than the ``true" SD which may be because of over-fitting. For the cv(Lasso) and cv(Lasso+OLS) adjusted estimator, the mean of Neyman SD estimator is within $1 \pm 7\%$ of the ``true" SD. Although the Neyman variance estimate is asymptotically conservative, the finite sample behavior of the Neyman SD estimate can be progressive for the regression-based adjusted estimator. However, if we increase the sample size $n$ to $1000$, almost all the confidence intervals are conservative.


We conduct more simulation examples to evaluate the conditions assumed for asymptotic normality of the Lasso adjusted estimator. We use the same simulation setup as above, but for simplicity, we generate the potential outcomes from linear model (set $\beta^{(a2)}=\beta^{(b2)}=0 $) and remove the effects of the hidden covariates $z_i$ in generating the error terms $\epsilon^{(a)}_i$ and $\epsilon^{(b)}_i$ and set $\rho=0, \ n_A=125$. We find that the distribution of the cv(Lasso) adjusted estimator may be non-normal when:

\begin{itemize}
\item[(1).] The covariates are generated from Gaussian distribution and the error terms do not satisfy second moment condition, e.g., being generated from $t$ distribution with one degree of freedom, see the upper two subplots of Figure~\ref{fig:cond-not-hold} (in the main text) for the histograms of unadjusted the cv(Lasso) adjusted estimators (the corresponding p-values of Kolmogorov--Smirnov testing for normality are less than $2.2e-16$).  \item[(2).] The covariates do not have bounded fourth moments, e.g., being generated from $t$ distribution with three degrees of freedom, see the lower two subplots of Figure~\ref{fig:cond-not-hold} (in the main text) for the histograms of unadjusted the cv(Lasso) adjusted estimators (again, the corresponding p-values of Kolmogorov--Smirnov testing for normality are less than $2.2e-16$).
\end{itemize}
These findings indicate that our moment condition (Condition \ref{cond:moment} and Remark 1) cannot be dramatically weakened. However, we also find that the cv(Lasso) adjusted estimator still has smaller SD and $\sqrt{\textnormal{MSE}}$ than the unadjusted estimator even when these moment conditions do not hold.


\section{The design matrix of the PAC data}
In the PAC data, there are 59 covariates (main effects) including 50 indicators which may be correlated with the outcomes. One of the main effects (called interactnew) has heavy tail, so we do the transform: $x \rightarrow \log ( |x|+1 )$ to make it look like normal distributed. We then centralize and standardize the non-indicator covariates. The quadratic terms (9 in total) of non-indicator covariates and two-way interactions between main effects (1711 in total) may also contribute to predict the potential outcomes, so we included them in the design matrix. The quadratic terms and the interactions between non-indicator covariates and the interactions between indicator covariates and non-indicator covariates are also centered and standardized. Some of the interactions are exactly the same as other effects and we only retain one of them. We also remove the interactions which are highly correlated (with correlation larger than $0.95$) with the main effects and remove the indicators with very sparse entries (where the number of 1's is less than 20). Finally, we form a design matrix $X$ with $1172$ columns (covariates) and $1013$ rows (subjects).

\section{Estimation of constants in the conditions}
Let $S^{(a)}=\{j: \bbeta^{(a)}_j \neq 0  \}$ and $S^{(b)}=\{j: \bbeta^{(b)}_j \neq 0  \}$ denote the sets of relevant covariates for treatment group and control group respectively. Denote $S=S^{(a)} \bigcup S^{(b)} = \{ j: \bbeta^{(a)}_j \neq 0 \ \textnormal{or} \  \bbeta^{(b)}_j \neq 0 \}$. We use bootstrap to get an estimation of the relevant covariates sets $S^{(a)}, S^{(b)}$ and then the approximation errors $e^{(a)}$ and $e^{(b)}$ are estimated by regressing the observed potential outcomes $a$ and $b$ on the covariates in $S$ respectively. We only present how to estimate $S^{(a)}$ and $e^{(a)}$ in detail and the estimation of $S^{(b)}$ and $e^{(b)}$ are similar. 

Let $A$, $B$ be the set of treated subjects (using PAC) and control subjects (without using PAC) respectively. Denote $a_i,i \in A$ the potential outcomes (quality-adjusted life years (QALYs)) under treatment and $x_i \in R^{1172}$ the covariates vector of the $i$th subject. For each $d=1,...,1000$, we draw a bootstrap sample $\{ (a_i^*(d),x_i^*(d)): i \in A \}$ with replacement from the data points $\{(a_i,x_i): i \in A\}$. Then computing the LassoOLS(CV) adjusted vector $\hat \bbeta(d)$ based on each bootstrap sample $\{ (a_i^*(d),x_i^*(d)): i \in A \}$. Let $\tau_j$ be the selection fraction of non-zero $\hat \beta_j(d)$ in the $1000$ bootstrap estimators, i.e., $\tau_j = (1/1000) \sum_{d=1}^{1000} \mathbb{I}_{ \{\hat \beta_j(d) \neq 0\} }$, where $\mathbb{I}$ is the indicator function. We form the relevant covariates $S^{(a)}$ by the covariates whose selection fraction are larger than $0.5$: $S^{(a)} = \{j:  \tau_j>0.5 \}$.

To estimate the approximation error $e^{(a)}$, we regress $a_i$ on the relevant covariates $x_{ij}, j \in S^{(a)}$ and compute OLS estimate and the corresponding residual. That is, let $T^{(a)}$ denote the complement set of $S^{(a)}$,
{\small \begin{equation}
 \bbeta^{(a)}_\textnormal{OLS} = \argmin_{\bbeta: \ \beta_j=0, \ \forall j \in  T^{(a)} } \frac{1}{2n_A}\sum_{i\in A} \left( a_i - \bar{a}_A  - ( \bx_i - \bar{\bx}_A )^T \bbeta \right)^2. \nonumber
\end{equation}}
\begin{equation}
 e^{(a)}_i = a_i - \bar{a}_A - ( \bx_i - \bar{\bx}_A )^T \bbeta^{(a)}_\textnormal{OLS}, \ i \in A.  \nonumber
\end{equation}

The maximal covariance $\delta_n$ is estimated as:
\begin{equation}
\begin{split}
 \max & \left\{ \max_j  \left| \frac{1}{n_A} \sum_{i \in A} \left( x_{ij}-(\bar{\bx})_j \right) \left( e^{(a)}_i - \bar{e}^{(a)}_A \right) \right|, \right. \\
& \left. \max_j  \left| \frac{1}{n_B} \sum_{i \in B} \left( x_{ij}-(\bar{\bx})_j \right) \left( e^{(b)}_i - \bar{e}^{(b)}_B \right) \right|
   \right\}.  \nonumber
\end{split}
\end{equation}

\section{Proofs of Theorem \ref{normality-thm}, \ref{conservative_variance}, \ref{conservative_variance-2} and Corollary \ref{compare-asym-var}}\label{A}
In this section, we will prove Theorem \ref{normality-thm} - \ref{conservative_variance-2} and Corollary \ref{compare-asym-var} under weaker sparsity conditions. 

\begin{defn}
We define an approximate sparsity measure. Given the regularization parameter $\lambda_a, \lambda_b$ and $\bbeta^{(a)}$ and $\bbeta^{(b)}$, the sparsity measures for treatment and control groups, $s^{(a)}_{\lambda_a}$ and $s^{(b)}_{\lambda_b}$ are defined as
\begin{equation}\label{def:s}
s^{(a)}_{\lambda_a} = \sum_{j=1}^p \min\left\{ \frac{ |\beta_j^{(a)}| }{\lambda_a}, 1 \right\}, \  s^{(b)}_{\lambda_b} = \sum_{j=1}^p \min\left\{ \frac{ |\beta_j^{(b)}| }{\lambda_b}, 1 \right\},
\end{equation}
respectively. We will allow $s^{(a)}_{\lambda_a}$ and $s^{(b)}_{\lambda_b}$ to grow with $n$, though the notation does not explicitly show this. Note that this is weaker than strict sparsity, as it allows $\bbeta^{(a)}$ and $\bbeta^{(b)}$ to
have 
many small non-zero entries.
\end{defn}

\noindent
{\bf Condition (*).} Suppose there exist $\bbeta^{(a)}$, $\bbeta^{(b)}$, $\lambda_a$ and $\lambda_b$ such that the conditions \ref{first-cond}, \ref{cond:moment}, \ref{cond:limit} and the following statements 1, 2, 3 hold simultaneously.
\begin{itemize}
\item {\bf Statement 1}. Decay and scaling. Let $s_{\lambda} = \max \left\{ s^{(a)}_{\lambda_a}, s^{(b)}_{\lambda_b} \right\}$,
\begin{equation}
\label{cond:delta_n-apx}
\delta_n = o\left( \frac{1}{s_{\lambda}\sqrt{\log p}} \right),
\end{equation}
\begin{equation}
\label{cond:s-scaling-apx}
(s_{\lambda} \log p)/{\sqrt n} = o(1).
\end{equation}

\item {\bf Statement 2}. Cone invertibility factor. Define the Gram matrix as
$\Sigma = n^{-1}\sum_{i=1}^n ( \bx_i - \bar{\bx} )( \bx_i - \bar{\bx} )^T$:
There exist constants $C >0$ and $\xi >1$ not depending on $n$, such that
\begin{equation}\label{cond:cone-inv}
\|\bh_S\|_1 \leq Cs_{\lambda} \| \Sigma \bh \|_\infty, \ \forall \bh \in \mathcal{C},
\end{equation}
with $\mathcal{C}=\{\bh: \|\bh_{S^c}\|_1 \leq \xi \|\bh_{S}\|_1\}$, and
\begin{equation} \label{def:S}
S = \{ j: |\beta^{(a)}_j| > \lambda_a \ \textnormal{or} \  |\beta^{(b)}_j| > \lambda_b \}.
\end{equation}

\item {\bf Statement 3}. Let $\tau = \min\big\{1/70,(3p_A)^2/70,(3-3p_A)^2/70 \big\}$. For constants $0< \eta < \frac{\xi - 1}{\xi + 1}$ and $0< M<\infty$, assume the regularization parameters of the Lasso belong to the sets
\begin{equation}\label{cond:lambda-a}
\lambda_a \in   (\frac{1}{\eta}, M] \times \left( \frac{2(1+\tau)L^{1/2}}{p_A} \sqrt{ \frac{2\log p}{n} } + \delta_n \right),
\end{equation}
\begin{equation} \label{cond:lambda-b}
\lambda_b \in   (\frac{1}{\eta}, M] \times \left( \frac{2(1+\tau)L^{1/2}}{p_B} \sqrt{ \frac{2\log p}{n} } + \delta_n \right).
\end{equation}

\end{itemize}

It is easy to verify that Condition (*) is implied by conditions \ref{first-cond} - \ref{last-last-cond}. In the following, we will prove Theorem \ref{normality-thm} - \ref{conservative_variance-2} and Corollary \ref{compare-asym-var} under the weaker Condition (*). For ease of notation, we will omit the subscript of $\hat \bbeta^{(a)}_\textnormal{Lasso}$, $\hat \bbeta^{(b)}_\textnormal{Lasso}$, $s_{\lambda}$, $s^{(a)}_{\lambda_a}$ and $s^{(b)}_{\lambda_b}$ from now on. Moreover, we can assume, without loss of generality, that
\begin{eqnarray}\label{zero-means}
\bar{a}=0, \ \bar{b}=0, \ \bar{\bx}=\mathbf{0}.
\end{eqnarray}
Otherwise, we can consider $\breve{a}_i = a_i - \bar{a}$, $\breve{b}_i = b_i - \bar{b}$ and $\breve{\bx}_i = \bx_i - \bar{\bx}$. Then, $\textnormal{ATE} = \bar{a} - \bar{b} = 0$ and the definition of $\ATElasso$ becomes
\begin{eqnarray}
\ATElasso  =  \left[ \bar{a}_A - ( \bar{\bx}_A )^T \hat{\bbeta}^{(a)} \right] -  \left[ \bar{b}_B - ( \bar{\bx}_B )^T \hat{\bbeta}^{(b)} \right].
\end{eqnarray}

We will rely heavily on the following Massart concentration inequality for sampling without replacement.

\begin{lem} \label{thm:concentration}
Let $\{z_i, i=1,...,n\}$ be a finite population of real numbers. Let $A \subset \{i,\ldots,n\}$ be a subset of deterministic size $|A|=n_A$ that is selected randomly without replacement. Define $p_A = n_A/n, \  \sigma^2 = \hbox{$n^{-1} \sum_{i=1}^{n}$} (z_i-\bar{z})^2$. Then, for any $t > 0$,
\begin{eqnarray}
\label{bernstein}
P\left(\bar{z}_A   - \bar{z} \geq t \right) \leq \exp \left\{ - \frac{p_A n_A t^2}{(1+\tau)^2 \sigma^2} \right\},
\end{eqnarray}
with $\tau = \min\big\{1/70,(3p_A)^2/70,(3-3p_A)^2/70 \big\}$.
\end{lem}
\noindent
{\bf Remark.} {\em Massart showed in his paper \cite{massart1986rates} that for sampling without replacement, the following concentration inequality holds:
\[ P\left(\bar{z}_A   - \bar{z} \geq t \right) \leq \exp \left\{ - \frac{p_A n_A t^2}{\sigma^2}  \right\}. \]
His proof required that $n/n_A$ must be an integer. We extend the proof to allow $n/n_A$ to be a non-integer but with a slightly larger constant factor $(1+\tau)^2$.
}

\begin{proof}
Assume $\bar{z}=0$ without loss of generality. For $n_A\le n/2$, let $m\ge 2$ and $r\ge 0$ be integers satisfying $n-n_Am = r < n_A$. Let $u\geq 0$. We first prove that
\begin{equation}\label{pf-lm-1-1}
\begin{split}
& E \exp\left( u \sum_{i\in A}z_i \right) \\
& \le E \exp\left( u \delta \sum_{i\in B}z_i/\{m(m+1)\} + u^2 n \sigma^2/4\right)
\end{split}
\end{equation}
for a random subset $B\subset\{1,\ldots,n\}$ of fixed size $|B|\le n/2$ and a certain fixed $\delta \in \{-1,1\}$. Let $P_1$ be the probability under which $\{i_1,\ldots,i_n\}$ is a random permutation of $\{1,\ldots,n\}$. Given $\{i_1,\ldots,i_n\}$, we divide the sequence into consecutive blocks $B_1,\ldots, B_{n_A}$ with $|B_j| = m+1$ for $j=1,\ldots, r$ and $|B_j| = m$ for $j=r+1,...,n_A$. Let $\bar{z}_k$ be the mean of $\{z_i: i \in B_k\}$ and $P_2$ be a probability conditionally on $\{i_1,\ldots,i_n\}$ under which $w_k$ is a random element of $\{z_i: i \in B_k\}$, $k=1,\ldots, n_A$. Then $\{w_1,\ldots,w_{n_A}\}$ is a random sample from $\{z_1,\ldots,z_n\}$ without replacement under $P=P_1P_2$. Let $\Delta_k = \max_{i\in B_k}z_i - \min_{i\in B_k}z_i$ and denote $E_2$ the expectation under $P_2$. The Hoeffding inequality gives
\begin{equation}
E_2 \exp\left( u \sum_{k=1}^{n_A}w_k\right) \le \exp\left( u \sum_{k=1}^{n_A}\bar{z}_k + (u^2/8)\sum_{k=1}^{n_A}\Delta_k^2 \right).
\end{equation}
As $\Delta_i^2 \le 2\sum_{i\in B_k}(z_i - \bar{z}_k)^2\le 2\sum_{i\in B_k}z_i^2$,
\begin{eqnarray}
E_2 \exp\left( u \sum_{k=1}^{n_A}w_k\right) \le \exp\left( u \sum_{k=1}^{n_A}\bar{z}_k + u^2 n \sigma^2/4\right)
\end{eqnarray}
Let $B = \cup_{k=1}^r B_k$. As $\bar{z} = 0$,
\begin{equation}
\sum_{k=1}^{n_A}\bar{z}_k = \sum_{i\in B}z_i/\{m(m+1)\}.
\end{equation}
This yields \eqref{pf-lm-1-1} with $\delta=1$ when $|B|\le n/2$. Otherwise, \eqref{pf-lm-1-1} holds with $\delta=-1$ when $B$ is replaced by $B^c$, as $\sum_{i\in B}z_i = - \sum_{i\in B^c}z_i$ due to $\bar{z}=0$.

Now, as $m(m+1)\ge 6$, repeated application of \eqref{pf-lm-1-1} yields
\begin{eqnarray}
 & & E \exp\left( u \sum_{i\in A}z_i \right) \nonumber \\
 &\le & E \exp\left[ u \delta' \sum_{i\in B'}z_i/\{m(m+1)m'(m'+1)\} \right. \nonumber \\
 & & \left. + \left(1+\{m(m+1)\}^{-2}\right)u^2 n \sigma^2/4 \right] \nonumber \\
&\le&  \exp\left[\left(1+\{m(m+1)\}^{-2}(1+1/36+1/36^2+ \right. \right. \nonumber \\
&& \left. \left. \cdots) \right)u^2 n \sigma^2/4\right] \nonumber \\
&=&  \exp\left[\left(1+(36/35)\{m(m+1)\}^{-2}\right)u^2 n \sigma^2/4\right] \nonumber \\
&\le &  \exp\left[\left(1+\tau\right)^2u^2 n \sigma^2/4\right]
\end{eqnarray}
with $\tau = (18/35)\{m(m+1)\}^{-2}$. The upper bound for $\tau$ follows from $2\le m < n/n_A< m+1$.

As $\bar{z}=0$, we also have
\begin{equation}
E \exp\left( u \sum_{i\in A}z_i \right) \le  \exp\left[\left(1+\tau \right)^2u^2 n \sigma^2/4\right]
\end{equation}
for $n_A>n/2$. This yields \eqref{bernstein} via the usual
\begin{eqnarray}
& &  P\left\{ \bar{z}_A - \bar{z} > t \right\} \nonumber \\
& \le &  \exp\left[ - u t +  (1+\tau)^2u^2 n \sigma^2/(4n_A^2)\right]  \nonumber \\
& = & \exp\left[ - 2 \frac{p_A n_A t^2}{(1+\tau)^2 \sigma^2} + \frac{p_A n_A t^2}{(1+\tau)^2 \sigma^2} \right]
\end{eqnarray}
with $u = 2p_An_At/\{\sigma(1+\tau)\}^2$.
\end{proof}

\subsection{Proof of Theorem \ref{normality-thm}}
\begin{proof}
Recall the decompositions of the potential outcomes:
\begin{equation}
a_i = \bar{a} + ( \bx_i - \bar{\bx} )^T \bbeta^{(a)} + \err^{(a)}_i = \bx_i^T \bbeta^{(a)} + \err^{(a)}_i,
\end{equation}
\begin{equation}
b_i = \bar{b} + ( \bx_i - \bar{\bx} )^T \bbeta^{(b)} + \err^{(b)}_i = \bx_i^T \bbeta^{(b)} + \err^{(b)}_i.
\end{equation}
If we define $\bh^{(a)} = \hat{\bbeta}^{(a)} - \bbeta^{(a)}$, $\bh^{(b)} = \hat{\bbeta}^{(b)} - \bbeta^{(b)}$, by substitution, we have
\begin{equation*}
\begin{split}
&\sqrt{n} ( \ATElasso - ATE ) \\
& =  \underbrace{ \sqrt{n} \left[ \bar{\err}_A^{(a)} - \bar{\err}_B^{(b)} \right] }_{\hypertarget{ATE-lhs}{}*} -  \underbrace{ \sqrt{n} \left[ \left( \bar{\bx}_A \right)^T{\bh^{(a)}} - \left( \bar{\bx}_B \right)^T{\bh^{(b)}} \right]}_{\hypertarget{ATE-rhs}{}**}.
\end{split}
\end{equation*}

We will analyze these two terms separately, showing that \hyperlink{ATE-lhs}{$(*)$} is asymptotically normal with mean $0$ and variance given by \eqref{sig-def}, and that \hyperlink{ATE-rhs}{$(**)$} is $o_p\left(1\right)$.

Asymptotic normality of \hyperlink{ATE-lhs}{$(*)$} follows from the Theorem 1 in \cite{Freedman2008b} with $a$ and $b$ replaced by $e^{(a)}$ and $e^{(b)}$ respectively. To bound \hyperlink{ATE-rhs}{$(**)$}, we will apply H\"older inequality to each of the terms. We will focus on the term involving the treatment group $A$, but exact same analysis is applied to the control group $B$. We have the bound
\begin{equation} \label{holder}
\left| \left( \bar{\bx}_A \right)^T{\bh^{(a)}} \right| \leq
\norm{ \bar{\bx}_A }_{\infty}
\shortnorm{ \bh^{(a)} }_1.
\end{equation}

We will bound the two terms on the right hand side of \eqref{holder} by the following Lemma \ref{lem:xterm} and Lemma \ref{lem:hterm}, respectively.
\begin{lem}
\label{lem:xterm}
Under the moment condition of [\ref{cond:xmoment}], if we let $c_n = \frac{(1+\tau)L^{1/4}}{p_A} \sqrt{\frac{2\log p}{n}}$, then as $n \rightarrow \infty$,
\begin{equation*}
P\left( \left\Vert \bar{\bx}_A \right\Vert_\infty > c_n \right) \rightarrow 0.
\end{equation*}
Thus, $\norm{ \bar{\bx}_A }_\infty = O_p\left( \sqrt{\frac{\log p}{n}}\right).$
\end{lem}

\begin{lem} \label{lem:hterm}
Assume the conditions of Theorem ~\ref{normality-thm} hold. Then $\shortnorm{ \bh^{(a)} }_1 =
o_p \left(\frac{1}{\sqrt{\log p}}\right) $.
\end{lem}

The proofs of these two Lemmas are below. 
Using these two Lemmas, it is easy to show that \hyperlink{ATE-rhs}{$(**)$}$=\sqrt{n}\cdot O_p \left( \sqrt{\frac{\log p}{n}}\right) \cdot o_p \left(\frac{1}{\sqrt{\log p}}\right)  = o_p\left(1\right)$.

\end{proof}

\subsection{Proof of Corollary \ref{compare-asym-var}}
\begin{proof}
By Theorem 1 in \cite{Freedman2008b}, the asymptotic variance of $\sqrt{n} \ \widehat{ATE}_\textnormal{unadj}$ is $\frac{1-p_A}{p_A} \lim_{n\rightarrow \infty} \sigma^2_{a} + \frac{p_A}{1-p_A} \lim_{n\rightarrow \infty} \sigma^2_{b} + 2 \lim_{n\rightarrow \infty} \sigma_{ab}$, so the difference is
\begin{align*}
& \frac{1-p_A}{p_A} \lim_{n\rightarrow \infty} \left( \sigma^2_{e^{(a)}} - \sigma^2_{a} \right) +
\frac{p_A}{1-p_A} \lim_{n\rightarrow \infty} \left( \sigma^2_{e^{(b)}} - \sigma^2_{b} \right) \\
& + 2 \lim_{n\rightarrow \infty} \left( \sigma_{e^{(a)}e^{(b)}} - \sigma_{ab} \right).
\end{align*}
We will analyze these three terms separately. Since $X\bbeta^{(a)}$ and $X\bbeta^{(b)}$ are the orthogonal projections of $a$ and $b$ onto the same subspace, we have
\[ (X\bbeta^{(a)})^Te^{(a)} = (X\bbeta^{(a)})^Te^{(b)} \]
\[ = (X\bbeta^{(b)})^Te^{(a)} = (X\bbeta^{(b)})^Te^{(b)} = 0. \]
Simple calculations yield
\begin{equation} \label{compare-asym-var-part1}
\sigma^2_{e^{(a)}} - \sigma^2_{a}  = ||e^{(a)}||_2^2 - ||a||_2^2 = - ||X\bbeta^{(a)}||_2^2, \nonumber
\end{equation}
\begin{equation} \label{compare-asym-var-part2}
\sigma^2_{e^{(b)}} - \sigma^2_{b}  = ||e^{(b)}||_2^2 - ||b||_2^2 = - ||X\bbeta^{(b)}||_2^2, \nonumber
\end{equation}
\begin{equation} \label{compare-asym-var-part3}
\sigma_{e^{(a)}e^{(b)}} - \sigma_{ab} = (e^{(a)})^T(e^{(b)}) - a^Tb = - (X\bbeta^{(a)})^T(X\bbeta^{(b)}). \nonumber
\end{equation}
Combining the above three equalities, we obtain the corollary.
\end{proof}

\subsection{Proof of Theorem \ref{conservative_variance}}
\begin{proof} To prove Theorem~\ref{conservative_variance}, it is enough to show that
\begin{equation}
\label{conver_vara}
\hat \sigma^2_{e^{(a)}} \stackrel{p}{\rightarrow}  \mathop {\lim}\limits_{n \rightarrow \infty }   \sigma^2_{e^{(a)}},
\end{equation}
\begin{equation}
\label{conver_varb}
\hat \sigma^2_{e^{(b)}} \stackrel{p}{\rightarrow}  \mathop {\lim}\limits_{n \rightarrow \infty } \sigma^2_{e^{(b)}}.
\end{equation}
We will only prove the statement \eqref{conver_vara} and omit the proof of the statement \eqref{conver_varb} since it is very similar.

We first state the following two lemmas. Lemma~\ref{lem:num-of-selcted-varibles} bounds the number of selected covariates (covariates with a nonzero coefficient), while Lemma~\ref{lem:subsamplemean} establishes conditions under which the subsample mean (without replacement) converges in probability to the population mean.
\begin{lem}
\label{lem:num-of-selcted-varibles}
Under the conditions in Theorem \ref{conservative_variance}, there exists a constant $C$, such that the following holds with probability going to 1:
\begin{equation}
 \hat s^{(a)} \leq C s; \ \ \hat s^{(b)} \leq C s.
\end{equation}
\end{lem}
The proof of Lemma~\ref{lem:num-of-selcted-varibles} can be found below. 

\begin{lem}
\label{lem:subsamplemean}
Let $\{z_i, i=1,...,n\}$ be a finite population of real numbers. Let $A \subset \{i,\ldots,n\}$ be a subset of deterministic size $ |A| = n_A$ that is selected randomly without replacement. Suppose that the population mean of the $z_i$ has a finite limit and that there exist constants $\epsilon>0$ and $L<\infty$ such that
\begin{equation}
\label{con:eps}
\frac{1}{n} \sum_{i=1}^{n} |z_i|^{1+\epsilon} \leq L.
\end{equation}
If $\frac{n_A}{n} \rightarrow p_A \in (0,1)$, then
\begin{equation}
\bar{z}_A  \stackrel{p}{\rightarrow}  \mathop {\lim}\limits_{n \rightarrow \infty } \bar{z}.
\end{equation}
\end{lem}
By definition \eqref{var_estim_a} and simple calculations,
\begin{eqnarray}
 && \hat \sigma^2_{e^{(a)}} \nonumber \\
 & = & \frac{1}{n_A - df^{(a)} } \sum_{i \in A} \left( a_i - \bar{a}_A - ( \bx_i - \bar{\bx}_A )^T \hat \bbeta^{(a)} \right)^2  \nonumber \\
 & = &  \frac{1}{n_A - df^{(a)} } \sum_{i \in A} \left( a_i - \bar{a}_A - ( \bx_i - \bar{\bx}_A )^T  \bbeta^{(a)} \right. \nonumber \\
 && \left. + ( \bx_i - \bar{\bx}_A )^T ( \bbeta^{(a)} - \hat \bbeta^{(a)}  ) \right)^2 \nonumber \\
 & = & \frac{1}{n_A - df^{(a)} } \sum_{i \in A} \left( a_i - \bx_i^T\bbeta^{(a)} - ( \bar{a}_A  - ( \bar{\bx}_A )^T  \bbeta^{(a)}) \right.  \nonumber \\
 & & \left. + ( \bx_i - \bar{\bx}_A )^T ( \bbeta^{(a)} - \hat \bbeta^{(a)}  ) \right)^2 \nonumber \\
 & = & \frac{n_A}{n_A - df^{(a)} }  \frac{1}{n_A} \sum_{i \in A} \left( e^{(a)}_i  -  \bar{e}_A^{(a)} + ( \bx_i - \bar{\bx}_A )^T \right.  \nonumber \\
 & & \left. ( \bbeta^{(a)} - \hat \bbeta^{(a)}  ) \right)^2 \nonumber \\
 & = & \frac{n_A}{n_A - df^{(a)} } \left\{  \frac{1}{n_A} \sum_{i \in A} \left( e^{(a)}_i  -  \bar{e}_A^{(a)} \right)^2 + \frac{1}{n_A} \sum_{i \in A}  \right. \nonumber \\
 & & \left. \left( ( \bx_i - \bar{\bx}_A )^T ( \bbeta^{(a)} - \hat \bbeta^{(a)}  ) \right)^2 \right\} + \frac{n_A}{n_A - df^{(a)} } \nonumber \\
 & & \left\{ \frac{1}{n_A} \sum_{i \in A}  ( e^{(a)}_i  -  \bar{e}_A^{(a)} )( \bx_i - \bar{\bx}_A )^T ( \bbeta^{(a)} - \hat \bbeta^{(a)}  ) \right\}. \nonumber
\end{eqnarray}
The second to last equality is due to the decomposition of potential outcome $a$:
\[  a_i = \bx_i^T \bbeta^{(a)} + e^{(a)}_i;  \ \  \bar{a}_A = ( \bar{\bx}_A )^T \bbeta^{(a)} +  \bar{e}_A^{(a)}. \]
It is easy to see that
\begin{equation}
\frac{1}{n_A} \sum_{i \in A} \left( e^{(a)}_i -  \bar{e}_A^{(a)} \right)^2 = \frac{1}{n_A} \sum_{i \in A} ( e^{(a)}_i)^2 - (\bar{e}_A^{(a)})^2.
\end{equation}
By the $4th$ moment condition on the approximation error $e^{(a)}$ (see \eqref{cond:errmoment}), and applying Lemma \ref{lem:subsamplemean} we get
\[  \frac{1}{n_A} \sum_{i \in A} ( e^{(a)}_i)^2 \stackrel{p}{\rightarrow} \mathop {\lim}\limits_{n \rightarrow \infty } \sigma^2_{e^{(a)}}; \ \ \bar{e}_A^{(a)} \stackrel{p}{\rightarrow} \mathop {\lim}\limits_{n \rightarrow \infty }  \bar e^{(a)} =0. \]
Therefore,
\begin{equation}
\label{part1}
\frac{1}{n_A} \sum_{i \in A} \left( e^{(a)}_i   -  \bar{e}_A^{(a)} \right)^2 \stackrel{p}{\rightarrow} \mathop {\lim}\limits_{n \rightarrow \infty } \sigma^2_{e^{(a)}}.
\end{equation}
Simple algebra operations give
\begin{eqnarray}
\label{secondterm}
& & \frac{1}{n_A} \sum_{i \in A} \left( (\bx_i - \bar{\bx}_A)^T ( \bbeta^{(a)} - \hat \bbeta^{(a)}  ) \right)^2  \nonumber \\
& = & ( \bbeta^{(a)} - \hat \bbeta^{(a)}  )^T \left[ \frac{1}{n_A} \sum_{i \in A} (\bx_i - \bar{\bx}_A)(\bx_i - \bar{\bx}_A)^T \right] \nonumber \\
& & ( \bbeta^{(a)} - \hat \bbeta^{(a)}  ) \nonumber \\
& \leq & || \bbeta^{(a)} - \hat \bbeta^{(a)} ||_1^2  \cdot || \frac{1}{n_A} \sum_{i \in A} (\bx_i - \bar{\bx}_A)(\bx_i - \bar{\bx}_A)^T  ||_\infty. \nonumber \\
\end{eqnarray}
We next show that \eqref{secondterm} converges to $0$ in probability. By Lemma~\ref{lem:hterm} and Lemma~\ref{lem:cov-mat-bound}, we have
\begin{equation} \label{secondterm:1}
|| \bbeta^{(a)} - \hat \bbeta^{(a)} ||_1 = \shortnorm{ \bh^{(a)} }_1 = o_p \left(\frac{1}{\sqrt{\log p}} \right),
\end{equation}
\begin{equation} \label{secondterm:2}
|| \frac{1}{n_A} \sum_{i \in A} (\bx_i - \bar{\bx}_A)(\bx_i - \bar{\bx}_A)^T  ||_\infty = O_p(1).
\end{equation}
Therefore,
\begin{equation}
\label{part2}
\frac{1}{n_A} \sum_{i \in A} \left( (\bx_i - \bar{\bx}_A)^T ( \bbeta^{(a)} - \hat \bbeta^{(a)}  ) \right)^2 \stackrel{p}{\rightarrow} 0.
\end{equation}
By Cauchy-Schwarz inequality,
\begin{eqnarray}
\label{part3}
& & | \frac{1}{n_A} \sum_{i \in A}  ( e^{(a)}_i  -  \bar{e}_A^{(a)} )(\bx_i - \bar{\bx}_A)^T ( \bbeta^{(a)} - \hat \bbeta^{(a)}  ) | \nonumber \\
& \leq &  \left[ \frac{1}{n_A} \sum_{i \in A}  \left( e^{(a)}_i  -  \bar{e}_A^{(a)} \right)^2 \right]^{\frac{1}{2}} \cdot \nonumber \\
& & \left[   \frac{1}{n_A} \sum_{i \in A} \left( (\bx_i - \bar{\bx}_A)^T ( \bbeta^{(a)} - \hat \bbeta^{(a)}  ) \right)^2 \right]^{\frac{1}{2}}
\end{eqnarray}
which converges to $0$ in probability because of \eqref{part1} and \eqref{part2}.

By Lemma~\ref{lem:num-of-selcted-varibles} and Condition \ref{cond:scaling}, we have
\begin{equation}
\label{part4}
\frac{n_A}{n_A - df^{(a)} } = \frac{n_A}{n_A - \hat s^{(a)} -1 } \stackrel{p}{\rightarrow} 1.
\end{equation}
Combining \eqref{part1}, \eqref{part2}, \eqref{part3} and \eqref{part4}, we conclude that
\[ \hat \sigma^2_{e^{(a)}}  \stackrel{p}{\rightarrow} \mathop {\lim}\limits_{n \rightarrow \infty }   \sigma^2_{e^{(a)}}.\]

The remaining part of the proof is to study the difference between the conservative variance estimate and the true asymptotic variance:
\begin{eqnarray}
   && \left( \frac{1}{p_A} \mathop {\lim}\limits_{n \rightarrow \infty } \sigma^2_{e^{(a)}} + \frac{1}{1-p_A} \mathop {\lim}\limits_{n \rightarrow \infty } \sigma^2_{e^{(b)}} \right) -  \left( \frac{1-p_A}{p_A}  \right.  \nonumber \\
   & & \left. \mathop {\lim}\limits_{n \rightarrow \infty } \sigma^2_{e^{(a)}} + \frac{p_A}{1-p_A} \mathop {\lim}\limits_{n \rightarrow \infty } \sigma^2_{e^{(b)}} + 2 \mathop {\lim}\limits_{n \rightarrow \infty } \sigma_{e^{(a)}e^{(b)}} \right)  \nonumber \\
 & = & \mathop {\lim}\limits_{n \rightarrow \infty } \sigma^2_{e^{(a)}} + \mathop {\lim}\limits_{n \rightarrow \infty } \sigma^2_{e^{(b)}}  - 2 \mathop {\lim}\limits_{n \rightarrow \infty } \sigma_{e^{(a)}e^{(b)}}   \nonumber \\
 & = &  \mathop {\lim}\limits_{n \rightarrow \infty } \sigma^2_{e^{(a)}-e^{(b)}} \nonumber \\
 & = &  \mathop {\lim}\limits_{n \rightarrow \infty } \frac{1}{n} \sum_{i=1}^{n} \left( a_i - b_i - \bx_i^T (\bbeta^{(a)} - \bbeta^{(b)}) \right)^2.
\end{eqnarray}

\end{proof}

\subsection{Proof of Theorem \ref{conservative_variance-2}}
\begin{proof}
By Lemma~\ref{lem:num-of-selcted-varibles}, $\max{ (\hat s^{(a)},\hat s^{(b)} )} = o_p(\min{ (n_A,n_B) })$. Therefore, $( \hat \sigma^2_{e^{(a)}}, \hat \sigma^2_{e^{(b)}} )$ and $( (\hat \sigma^*)^2_{e^{(a)}}, (\hat \sigma^*)^2_{e^{(b)}} )$ have the same limits. The conclusion follows from Theorem~\ref{conservative_variance}.
\end{proof}

\section{Proofs of Lemmas}
\label{lem-proofs}

In this section, we will drop the superscript on $\bh$, $e$ and $\hat \bbeta$ and focus on the proof for treatment group A, as the same analysis can be applied to control group B.

\subsection{Proof of Lemma~\ref{lem:xterm}}
\begin{proof}
Let $c_n = \frac{(1+\tau)L^{1/4}}{p_A} \sqrt{\frac{2\log p}{n}}$. By the union bound,
\begin{equation}\label{eq:unionx}
\begin{split}
P\left( \left\Vert \bar{\bx}_A \right\Vert_\infty > c_n \right) & =
P\left( \max_{j=1,\ldots,p} \left| \frac{1}{n_A}\sum_{i \in A}x_{ij} \right| > c_n \right) \\
& \leq \sum_{j=1}^{p} P\left( \left| \frac{1}{n_A}\sum_{i \in A}x_{ij} \right| > c_n \right).
\end{split}
\end{equation}
By Cauchy-Schwarz inequality, we have
\begin{equation}\label{var-x}
 \frac{1}{n} \sum_{i=1}^{n} x_{ij}^2  \leq \left( \frac{1}{n} \sum_{i=1}^{n} x_{ij}^4  \right)^{\frac{1}{2}}  \left( \frac{1}{n} \sum_{i=1}^{n} 1^2  \right)^{\frac{1}{2}}  \leq \sqrt{L}.
\end{equation}
Substituting the concentration inequality \eqref{bernstein} into \eqref{eq:unionx},
\begin{equation*}
\begin{split}
P\left( \left\Vert \bar{\bx}_A \right\Vert_\infty > c_n \right) & \leq
2\exp \left\{ \log p - \frac{p_A n_A c_n^2 }{ (1+\tau)^2 L^{1/2} }  \right\}  \\
& = 2 \exp \left\{ - \log p \right\} \rightarrow 0.
\end{split}
\end{equation*}
\end{proof}

\subsection{Proof of Lemma~\ref{lem:hterm}}
\begin{proof} We start with the KKT condition, which characterizes the solution to the Lasso. Recall the definition of the Lasso estimator $\hat \bbeta$:
\[ \hat \bbeta =  \argmin_{\bbeta} \frac{1}{2n_A} \sum_{i \in A} \left( a_i - \bar{a}_A - (\bx_i -  \bar{\bx}_A )^T \bbeta \right)^2 + \lambda_a \norm{\bbeta}_1.  \]
The KKT condition for $\hat \bbeta$ is
\begin{equation}
\label{KKT}
 \frac{1}{n_A} \sum_{i\in A} (\bx_i -  \bar{\bx}_A) \left( a_i -  \bar{a}_A  - (\bx_i -  \bar{\bx}_A)^T \hat \bbeta \right)  = \lambda_a \mathbf{\kappa},
\end{equation}
where $\kappa$ is the subgradient of $||\bbeta||_1$ taking value at $\bbeta = \hat \bbeta$, i.e.,
\begin{equation}\label{eqn:KKT}
\begin{split}
\mathbf{\kappa} \in \partial || \bbeta ||_1 \left |_{ \bbeta = \hat \bbeta } \right.  \quad \textnormal{with} \quad
\left\{
\begin{aligned}
\kappa_j &\in [-1,1] \textnormal{ for } j \st \hat\beta_j = 0 \\
\kappa_j &= \textnormal{sign}(\hat\beta_j) \textnormal{ otherwise}
\end{aligned}
\right.
\end{split}
\end{equation}
Substituting $a_i$ by the decomposition \eqref{decom-a}, \eqref{KKT} becomes
\begin{equation}
\begin{split}
\label{KKT1}
 &\frac{1}{n_A} \sum_{i\in A} (\bx_i - \bar{\bx}_A) (\bx_i - \bar{\bx}_A)^T ( \bbeta - \hat \bbeta ) \\
 & + \frac{1}{n_A} \sum_{i\in A} (\bx_i - \bar{\bx}_A)( e_i - \bar{e}_A )   = \lambda_a \mathbf{\kappa}.
\end{split}
\end{equation}
Multiplying both sides of \eqref{KKT1} by $ - \bh^T = (\bbeta -  \hat \bbeta  )^T$, we have
\begin{eqnarray}
\label{basic-inequality}
 & & \frac{1}{n_A} \sum_{i\in A} \left( (\bx_i - \bar{\bx}_A)^T \bh \right)^2 -  \bh^T \frac{1}{n_A} \sum_{i\in A} (\bx_i - \bar{\bx}_A)( e_i - \bar{e}_A ) \nonumber \\
 & & = \lambda_a (\bbeta -  \hat \bbeta  )^T \mathbf{\kappa}  \leq  \lambda_a \left( \norm{\bbeta}_1 - \shortnorm{\hat\bbeta}_1\right) \nonumber
\end{eqnarray}
where the last inequality holds because
\[ \bbeta ^T \mathbf{\kappa} \leq || \bbeta ||_1 || \mathbf{\kappa} ||_\infty \leq || \bbeta ||_1 \ \  \textnormal{and} \ \ \hat \bbeta^T \mathbf{\kappa} = || \hat \bbeta ||_1. \]
Rearranging, and applying H\"older's inequality, we have
\begin{align*}
& \frac{1}{n_A} \sum_{i\in A} \left( (\bx_i - \bar{\bx}_A)^T \bh \right)^2  \\
&\leq \lambda_a \left( \norm{\bbeta}_1 - \shortnorm{\hat\bbeta}_1\right)
+ \bh^T \frac{1}{n_A} \sum_{i\in A} (\bx_i - \bar{\bx}_A)( e_i - \bar{e}_A ) \\
&\leq \lambda_a \left( \norm{\bbeta}_1 - \shortnorm{\hat\bbeta}_1\right)
+ \norm{\bh}_1
\underbrace{\norm{ \frac{1}{n_A} \sum_{i\in A} (\bx_i - \bar{\bx}_A)( e_i - \bar{e}_A ) }_\infty}_{{\hypertarget{term:xe}{}*}}
\end{align*}
To control the term \hyperlink{term:xe}{$(*)$}, we define the event $\mathcal{L} = \left\{ \hyperlink{term:xe}\textnormal{*} \leq \eta \lambda_a  \right\}$. The following Lemma \ref{lem:xe-term} shows that, with $\lambda_a$ defined appropriately, $\mathcal{L}$ holds with probability approaching 1. We will prove this Lemma later.
\begin{lem}\label{lem:xe-term}
Define

$\mathcal{L} = \left\{ \left\Vert \frac{1}{n_A} \sum_{i \in A} ( \bx_i - \bar{\bx}_A ) ( e_i - \bar{e}_A ) \right\Vert_\infty \leq \eta \lambda_a  \right\}$.

\noindent Then under the conditions of Theorem \ref{normality-thm},
$
P(\mathcal{L}) \rightarrow 1.
$
\end{lem}
On $\mathcal{L}$
\begin{align} \label{on-event-L}
\frac{1}{n_A} \sum_{i\in A} \left( (\bx_i - \bar{\bx}_A)^T \bh \right)^2  &\leq
\lambda_a \left( \norm{\bbeta}_1 - \shortnorm{\hat\bbeta}_1 + \eta \norm{\bh}_1 \right).
\end{align}
By substituting the defiition of $\bh$, and several applications of the triangle inequality, we have
\begin{align*}
\norm{\bbeta}_1 - \shortnorm{\hat\bbeta}_1 \leq \norm{\bh_S}_1 - \norm{\bh_{S^c}}_1 + 2\norm{\bbeta_{S^c}}_1.
\end{align*}
Therefore,
\begin{align*}
& \frac{1}{n_A} \sum_{i\in A} \left( (\bx_i - \bar{\bx}_A)^T \bh \right)^2 \\
&\leq \lambda_a \left( \norm{\bh_S}_1 - \norm{\bh_{S^c}}_1 + 2\norm{\bbeta_{S^c}}_1 + \eta \norm{\bh}_1 \right) \\
&\leq \lambda_a \left( ( \eta - 1 ) \norm{\bh_{S^c}}_1 + (1 + \eta)\norm{\bh_S}_1 + 2 \norm{\bbeta_{S^c}}_1 \right).
\end{align*}
Because $\frac{1}{n_A} \sum_{i\in A} \left( (\bx_i - \bar{\bx}_A)^T \bh \right)^2 \geq 0$, we obtain
\begin{equation}\label{eqn:hscbound}
\begin{split}
& ( 1 - \eta ) \norm{\bh_{S^c}}_1 \\
& \leq (1 + \eta)\norm{\bh_S}_1 + 2 \norm{\bbeta_{S^c}}_1 \leq (1 + \eta)\norm{\bh_S}_1 + 2 s \lambda_a.
\end{split}
\end{equation}
where the last inequality holds because of the definition of $s$ in \eqref{def:s} and $S$ in \eqref{def:S}.

Consider the following two cases:

(I) If $(1+\eta)\shortnorm{\bh_S}_1 + 2s\lambda_a \geq (1-\eta)\xi\shortnorm{\bh_S}_1$ then by \eqref{eqn:hscbound},
\begin{equation*}
\begin{split}
\norm{\bh}_1  & = \norm{\bh_S}_1 + \norm{\bh_{S^c}}_1 \\
& \leq \left( \frac{1+\eta}{1-\eta} + 1\right)\shortnorm{\bh_S}_1 + \frac{2s\lambda_a}{1-\eta} \\
& \leq \frac{2s\lambda_a}{1-\eta}\left( \frac{2}{(1-\eta)\xi - (1+\eta)} + 1 \right).
\end{split}
\end{equation*}
By the definition of $\lambda_a$ and the scaling assumptions \eqref{cond:delta_n-apx}, \eqref{cond:s-scaling-apx}, we have that $s \lambda_a = o\left( \frac{1}{\sqrt{\log p}}\right)$.

(II) If $(1+\eta)\shortnorm{\bh_S}_1 + 2s\lambda_a < (1-\eta)\xi\shortnorm{\bh_S}_1$ then by  \eqref{eqn:hscbound} we have
$
\shortnorm{\bh_{S^c}}_1 \leq \xi \shortnorm{\bh_S}_1
$.
Applying the cone invertibility condition on the design matrix \eqref{cond:cone-inv},
\begin{equation}\label{eqn:hbound}
\begin{split}
\norm{\bh}_1 & = \norm{\bh_S}_1 + \norm{\bh_{S^c}}_1 \\
& \leq (1+\xi)\shortnorm{\bh_S}_1 \leq (1+\xi) C s \norm{ \frac{1}{n}X^TX \bh}_\infty
\end{split}
\end{equation}
Before applying this inequality we will revisit the KKT condition \eqref{eqn:KKT}, but this time we will take the $l_\infty$-norm, yielding
\begin{equation}\label{eqn:KKT-infty}
\begin{split}
&\norm{ \frac{1}{n_A} \sum_{i\in A} (\bx_i - \bar{\bx}_A) (\bx_i - \bar{\bx}_A)^T \bh }_\infty \\
&\leq \lambda_a + \norm{ \frac{1}{n_A} \sum_{i\in A} (\bx_i - \bar{\bx}_A)( e_i - \bar{e}_A ) }_\infty
\leq (1+\eta) \lambda_a,
\end{split}
\end{equation}
where the latter inequality holds on the set $\mathcal{L}$.
The final step is to control the deviation of the subsampled covariance matrix from the population covariance matrix, so that we can apply \eqref{eqn:hbound}. We define another event with constant $C_1 = \frac{2(1+\tau)L^{1/2}}{p_A} $
\begin{equation*}
\begin{split}
\mathcal{M} = & \left\{
\norm{  \frac{1}{n_A} \sum_{i\in A} (\bx_i - \bar{\bx}_A) (\bx_i - \bar{\bx}_A)^T   - \frac{1}{n}X^TX }_\infty \right. \\
& \left. \leq C_1  \sqrt{ \frac{\log p}{n} }
\right\}
\end{split}
\end{equation*}

\begin{lem}\label{lem:cov-mat-bound}
Assume stability of treatment assignment probability condition \ref{first-cond} and moment condition \ref{cond:xmoment} hold. Then $P(\mathcal{M}) \rightarrow 1$.
\end{lem}

We will prove Lemma \ref{lem:cov-mat-bound} later. Continuing our inequalities, on the event $\mathcal{L} \cap \mathcal{M}$,
\begin{equation*}
\begin{split}
& s \norm{ \frac{1}{n}X^TX \bh}_\infty \\
\leq & C_1  s \sqrt{ \frac{\log p}{n} } \norm{\bh}_1 + s \norm{\frac{1}{n_A} \sum_{i\in A} (\bx_i - \bar{\bx}_A) (\bx_i - \bar{\bx}_A)^T \bh}_\infty \\
\leq & o(1) \norm{\bh}_1 + s(1+\eta)\lambda_a,
\end{split}
\end{equation*}
where we have applied the scaling assumption \eqref{cond:s-scaling-apx} and \eqref{eqn:KKT-infty} in the second line. Hence, by \eqref{eqn:hbound},
\begin{equation*}
\norm{\bh}_1 \leq (1+\xi)C\left[ o(1)\norm{\bh}_1 + s(1+\eta)\lambda_a \right].
\end{equation*}
Again, applying the scaling assumptions \eqref{cond:delta_n-apx} and \eqref{cond:s-scaling-apx}, we get $\norm{\bh}_1 = o_p\left( \frac{1}{\sqrt{\log p}}\right)$.
\end{proof}

\subsection{Proof of Lemma~\ref{lem:num-of-selcted-varibles}}
\begin{proof}
In the proof of Lemma~\ref{lem:hterm}, we have shown that, on $\mathcal{L}$ defined in Lemma~\ref{lem:xe-term},
\begin{align} \label{on-event-L}
& \frac{1}{n_A} \sum_{i\in A} \left( (\bx_i - \bar{\bx}_A)^T ( \bbeta - \hat \bbeta ) \right)^2 \\
\leq & \lambda_a \left( \norm{\bbeta}_1 - \shortnorm{\hat\bbeta}_1 + \eta ||  \bbeta - \hat \bbeta  ||_1 \right). \nonumber \\
\leq & \lambda_a ( 1 + \eta  ) || \bbeta - \hat \bbeta ||_1.
\end{align}
Let $\bx^j$ be the $j$-th column of the design matrix $X$ and $\bar{\bx}_A^j = n_A^{-1} \sum_{i\in A} x_{ij}$. Again, by KKT conditon, we have
\[  \left| \frac{1}{n_A} \sum_{i\in A} (x_{ij} - \bar{\bx}_A^j )  \left( a_i - \bar{a}_A  - (\bx_i - \bar{\bx}_A)^T \hat \bbeta \right)  \right| = \lambda_a, \]
\[ \textnormal{if} \ \hat \bbeta_j \neq 0.   \]
Substituting $a_i$ by the decomposition \eqref{decom-a} yields
\begin{equation*}
\begin{split}
 & \left| \frac{1}{n_A} \sum_{i\in A} (x_{ij} - \bar{\bx}_A^j ) ( e_i - \bar{ e}_A ) + \frac{1}{n_A} \sum_{i\in A} (x_{ij} - \bar{\bx}_A^j ) \right. \\
 & \left. (\bx_i - \bar{\bx}_A)^T ( \bbeta - \hat \bbeta ) \right| = \lambda_a.
\end{split}
\end{equation*}
Combining with the definition of the event $\mathcal{L}$, we have if $\hat \bbeta_j \neq 0$
\begin{equation}
\label{lowerbound}
 \Delta_j := \left| \frac{1}{n_A} \sum_{i\in A} (x_{ij} - \bar{\bx}_A^j )  (\bx_i - \bar{\bx}_A)^T ( \bbeta - \hat \bbeta ) \right|  \geq (1-\eta) \lambda_a.
\end{equation}
Let $Z=(\bz_1,...,\bz_n) \in R^{p\times n}$ with $\bz_i = \bx_i - \bar{\bx}_A \in R^p$ and denote $\bw = Z^T ( \bbeta - \hat \bbeta )$, then
\begin{equation*}
\begin{split}
 \frac{1}{n_A} || \bw_A ||_2^2 & =  \frac{1}{n_A} \sum_{i\in A} \left( (\bx_i - \bar{\bx}_A)^T ( \bbeta - \hat \bbeta ) \right)^2 \\
 & \leq \lambda_a ( 1 + \eta ) || \bbeta - \hat \bbeta ||_1.
\end{split}
\end{equation*}
Let $Z_A = (\bz_i: i \in A)$; since the largest eigenvalues of $Z_A^T Z_A$ and $Z_A Z_A^T$ are the same,
\begin{eqnarray}
\label{upperbound}
 & & \frac{1}{n_A^2} \bw_A^T Z_A^T Z_A \bw_A \nonumber \\
 & \leq & \frac{1}{n_A^2} \lambda_\textnormal{max}( Z_A^T Z_A ) ||\bw_A||_2^2 \nonumber \\
 & \leq &  \frac{1}{n_A} \lambda_\textnormal{max}( Z_A Z_A^T ) \lambda_a ( \eta + 1 ) || \bbeta - \hat \bbeta ||_1 \nonumber \\
 & \leq & \Lambda_\textnormal{max} \frac{n}{n_A} \lambda_a ( 1 + \eta ) || \bbeta - \hat \bbeta ||_1. \nonumber
\end{eqnarray}
The last inequality holds because
\begin{eqnarray}
 & & \lambda_\text{max}( Z_A Z_A^T ) \nonumber \\
 & = & \max_{\bu: ||\bu||_2=1} \bu^T Z_A Z_A^T \bu \nonumber \\
 & = &  \max_{\bu: ||\bu||_2=1} \bu^T \sum_{i\in A} (\bx_i - \bar{\bx}_A) (\bx_i - \bar{\bx}_A)^T \bu \nonumber \\
 & = & \max_{\bu: ||\bu||_2=1} \bu^T \sum_{i\in A} \bx_i \bx_i^T \bu - n_A \bu^T (\bar{\bx}_A)(\bar{\bx}_A)^T \bu  \nonumber \\
 & \leq & \max_{\bu: ||\bu||_2=1} \bu^T \sum_{i\in A} \bx_i \bx_i^T \bu \leq n\Lambda_\textnormal{max} .
\end{eqnarray}
On the other hand,
\begin{equation}
\label{upperbound}
 \frac{1}{n_A^2} \bw_A^T Z_A^T Z_A \bw_A =  \sum_{j=1}^{p} \Delta_j^2 \geq   \sum_{j: \hat \beta_j \neq 0} \Delta_j^2 \geq  (1-\eta)^2\lambda_a^2 \hat s.
\end{equation}
Combining \eqref{lowerbound}, \eqref{upperbound} and the fact that with probability going to $1$ (see the proof of Lemma \ref{lem:hterm})
\[  || \bbeta - \hat \bbeta ||_1 \leq C s ( 1 + \eta ) \lambda_a, \]
where $C$ is a constant, we conclude that with probability going to $1$,
\begin{equation*}
\begin{split}
 \hat s & \leq \frac{1}{(1-\eta)^2} \frac{1}{\lambda_a^2}\Lambda_\textnormal{max} \frac{n}{n_A} \lambda_a ( 1 + \eta ) C s ( 1 + \eta ) \lambda_a \\
 & \leq \frac{C( 1+\eta )^2}{p_A(1-\eta)^2} s.
\end{split}
\end{equation*}

\end{proof}

\subsection{Proof of Lemma~\ref{lem:subsamplemean}}
\begin{proof}
For any $t>0$, we have
\begin{equation}
\label{sum}
 P( | \bar{z}_A - \mathop {\lim}\limits_{n \rightarrow \infty }  \bar{z} |  > t ) \leq P( | \bar{z}_A -  \bar{z} | > t/2 )  + P( | \bar{z} - \mathop {\lim}\limits_{n \rightarrow \infty } \bar{z} |  > t/2 ).
\end{equation}

The second term in the right hand side of \eqref{sum} obviously converges to $0$ as $n\rightarrow \infty$. To bound the first term, we apply the concentration inequality \eqref{bernstein}.
By \eqref{con:eps}, it is easy to show
\begin{equation*}
\begin{split}
 \frac{1}{n} \sum_{i=1}^{n} z_i^2 & = \frac{1}{n} \sum_{i=1}^{n} |z_i|^{1-\epsilon} |z_i|^{1+\epsilon}\\
  & \leq (nL)^{\frac{1-\epsilon}{1+\epsilon}} \frac{1}{n} \sum_{i=1}^{n} |z_i|^{1+\epsilon} \leq L^{\frac{2}{1+\epsilon}}n^{\frac{1-\epsilon}{1+\epsilon}}.
\end{split}
\end{equation*}
Concentration inequality \eqref{bernstein} yields
\[  P( | \bar{z}_A -  \bar{z} |  > t/2 ) \leq  2 \exp \left\{ - \frac{p_An_A t^2}{4(1+\tau)^2L^{\frac{2}{1+\epsilon}}n^{\frac{1-\epsilon}{1+\epsilon}} }    \right\}  \rightarrow 0. \]

\end{proof}

\subsection{Proof of Lemma~\ref{lem:xe-term}}

\begin{proof}
It is easy to verify that
\[ \frac{1}{n_A} \sum_{i \in A} ( \bx_i - \bar{\bx}_A ) ( e_i - \bar{e}_A ) = \frac{1}{n_A} \sum_{i \in A} \bx_i e_i -  (\bar{\bx}_A) (\bar{e}_A). \]
Hence,
\begin{equation} \label{left-side}
\begin{split}
& || \frac{1}{n_A} \sum_{i \in A} ( \bx_i - \bar{\bx}_A ) ( e_i - \bar{e}_A ) ||_\infty \\
& \leq || \frac{1}{n_A} \sum_{i \in A} \bx_i e_i ||_\infty + || (\bar{\bx}_A) (\bar{e}_A) ||_\infty.
\end{split}
\end{equation}
We analyze these two terms on the right hand side of the inequality separately. For the first term, by the triangle inequality and the definition of $\delta_n$ in \eqref{def:delta},
\begin{align} \label{xterm:part1}
& || \frac{1}{n_A} \sum_{i \in A} \bx_i e_i ||_\infty \\
& \leq || \frac{1}{n_A} \sum_{i \in A} \bx_i e_i -  \frac{1}{n} \sum_{i=1}^{n} \bx_i e_i ||_\infty + ||  \frac{1}{n} \sum_{i=1}^{n} \bx_i e_i ||_\infty  \nonumber \\
&\leq || \frac{1}{n_A} \sum_{i \in A} \bx_i e_i -  \frac{1}{n} \sum_{i=1}^{n} \bx_i e_i ||_\infty + \delta_n.
\end{align}
We will again bound \eqref{xterm:part1} by the concentration inequality \eqref{bernstein} in Lemma~\ref{thm:concentration}. By the Cauchy-Schwarz inequality, we have for any $j=1,..,p$,
\[ \frac{1}{n} \sum_{i=1}^{n} x_{ij}^2 e_i^2 \leq \left( \frac{1}{n} \sum_{i=1}^{n} x_{ij}^4 \right)^{\frac{1}{2}} \left( \frac{1}{n} \sum_{i=1}^{n} e_i^4 \right)^{\frac{1}{2}} \leq L.  \]
Let $t_n = \frac{(1+\tau)L^{1/2}}{p_A} \sqrt{ \frac{2\log p}{n} }$, then by the union bound and the concentration inequality \eqref{bernstein},
\begin{align*}
 & P\left( || \frac{1}{n_A} \sum_{i \in A} \bx_i e_i -  \frac{1}{n} \sum_{i=1}^{n} \bx_i e_i ||_\infty > t_n \right) \\
& \leq 2\exp \left\{ \log p - \frac{p_A n_A t_n^2}{(1+\tau)^2L)}  \right\} \\
& = 2 \exp \left\{ - \log p \right\} \rightarrow 0.
\end{align*}
Taking this back to \eqref{xterm:part1}, we have
\begin{equation} \label{xterm:part1-final}
P\left( || \frac{1}{n_A} \sum_{i \in A} \bx_i e_i ||_\infty  \leq  t_n + \delta_n \right)  \rightarrow 1.
\end{equation}

For the second term, by Lemma \ref{lem:xterm}, we have shown that,
\[ P\left( \left\Vert \bar{\bx}_A \right\Vert_\infty \leq  \frac{(1+\tau)L^{1/4}}{p_A} \sqrt{\frac{2\log p}{n}} \right) \rightarrow 1.  \]
A similar proof yields
\[ P\left( \left\Vert \bar{e}_A \right\Vert_\infty \leq  \frac{(1+\tau)L^{1/4}}{p_A} \sqrt{\frac{2\log p}{n}} \right) \rightarrow 1.  \]
Hence, under the scaling condition \eqref{cond:s-scaling-apx},
\begin{equation} \label{xterm:part2}
P\left( \left\Vert (\bar{\bx}_A) (\bar{e}_A) \right\Vert_\infty \leq  \frac{(1+\tau)L^{1/2}}{p_A} \sqrt{\frac{2\log p}{n}} \right) \rightarrow 1.
\end{equation}

Combining \eqref{xterm:part1-final} and \eqref{xterm:part2} yields
\begin{align*}
 P & \left(  || \frac{1}{n_A} \sum_{i \in A} ( \bx_i - \bar{\bx}_A ) ( e_i - \bar{e}_A ) ||_\infty \right.  \\
 & \left. \leq  \frac{2(1+\tau)L^{1/2}}{p_A} \sqrt{ \frac{2\log p}{n} } + \delta_n \right) \rightarrow 1.
\end{align*}
The conclusion follows from the condition $\lambda_a \in  (\frac{1}{\eta}, M] \times \left( \frac{2(1+\tau)L^{1/2}}{p_A} \sqrt{ \frac{2\log p}{n} } + \delta_n\right)$.
\end{proof}

\subsection{Proof of Lemma~\ref{lem:cov-mat-bound}}
\begin{proof}
It is easy to see that
\[ \frac{1}{n_A} \sum_{i\in A} (\bx_i - \bar{\bx}_A) (\bx_i - \bar{\bx}_A)^T = \frac{1}{n_A} \sum_{i\in A} \bx_i \bx_i^T - (\bar{\bx}_A)(\bar{\bx}_A)^T. \]
Then, by triangle inequality,
\begin{align}
& || \frac{1}{n_A} \sum_{i\in A} (\bx_i - \bar{\bx}_A) (\bx_i - \bar{\bx}_A)^T   - \frac{1}{n}X^TX ||_\infty  \\
& \underbrace{ \leq || \frac{1}{n_A} \sum_{i\in A} \bx_i \bx_i^T   - \frac{1}{n}\sum_{i=1}^{n} \bx_i \bx_i^T  ||_\infty }_{\hypertarget{mat-lhs}{}*} + \underbrace{ || (\bar{\bx}_A)(\bar{\bx}_A)^T ||_\infty }_{\hypertarget{mat-rhs}{}**}.
\end{align}

We control the first term \hyperlink{mat-lhs}{$(*)$} again using the concentration inequality \eqref{bernstein} and the union bound. By the wayCauchy-Schwarz inequality, for $j,k = 1,...,p$,
\[ \frac{1}{n}\sum_{i=1}^{n} x_{ij}^2x_{ik}^2 \leq \left( \frac{1}{n}\sum_{i=1}^{n} x_{ij}^4 \right)^{ \frac{1}{2} }  \left( \frac{1}{n}\sum_{i=1}^{n} x_{ik}^4 \right)^{ \frac{1}{2} } \leq L.  \]
Then,
\begin{eqnarray} \label{cov-mat-part1}
 && P\left(  || \frac{1}{n_A} \sum_{i\in A} \bx_i \bx_i^T   - \frac{1}{n}\sum_{i=1}^{n} \bx_i \bx_i^T  ||_\infty \geq \frac{(1+\tau)L^{1/2}}{p_A} \right. \nonumber \\
 && \left. \ \ \ \ \  \sqrt{\frac{3\log p}{n}} \right)  \nonumber \\
& \leq & 2 \exp\left\{ 2\log p - \frac{3 p_A n_A (1+\tau)^2L \log p }{(1+\tau)^2L p_A^2 n}   \right\} \nonumber \\
&  = & 2 \exp\left\{  - \log p \right\} \rightarrow 0.
\end{eqnarray}

The second term \hyperlink{mat-lhs}{$(**)$} is bounded by again observing that, by Lemma \ref{lem:xterm} and the scaling condition \eqref{cond:s-scaling-apx},
\begin{equation} \label{cov-mat-part2}
\hyperlink{mat-lhs}{(**)} \leq || \bar{\bx}_A ||_\infty^2 = o_p( \sqrt{ \frac{\log p}{n} } ).
\end{equation}
Combining \eqref{cov-mat-part1} and \eqref{cov-mat-part2} yields the conclusion.
\end{proof}

\section{Tables and Figures}

\begin{algorithm*}[htb]
\caption{\hspace{0.2cm} $K$-fold Cross Validation (CV) for the Lasso+OLS estimator}
\label{alg:cv}
\begin{algorithmic}[1]\vspace{0.2cm}
\REQUIRE
    Design matrix $X$, response $Y$ and a sequence of tuning parameter $\lambda_{1},...,\lambda_{J}$; Number of folds $K$.\vspace{0.2cm}
\ENSURE
    The optimal tuning parameter selected by CV: $\lambda_{optimal}$. \vspace{0.2cm}
\STATE Divide randomly the data $z=(X,Y)$ into $K$ roughly equal parts $z_k,k=1,...,K$; \label{code:fram:extract}
\STATE For each $k=1,...,K$, denote $\hat S^{(k)}(\lambda_{0}) = \emptyset $ and $\hat \beta^{(k)}_\textnormal{Lasso+OLS}(\lambda_{0}) = 0$.
\begin{itemize}
\item Fit the model with parameters $\lambda_{j},j=1,...,J$ to the other $K-1$ parts $z_{-k} = z\setminus z_k$ of the data, giving the Lasso solution path $\hat \beta^{(k)}(\lambda_{j}),j=1,...,J$ and compute the selected covariates set $\hat S^{(k)}(\lambda_{j}) = \{l:  \hat \beta^{(k)}_l(\lambda_{j}) \neq 0 \}, j=1,...,J$ on the path;
\item For each $j=1,...,J$, compute the Lasso+OLS estimator:
\begin{equation}
\hat \beta^{(k)}_\textnormal{Lasso+OLS}(\lambda_{j}) = \left\{
\begin{aligned}
& \argmin_{\bbeta: \ \beta_j=0, \ \forall j \notin \hat S^{(k)}(\lambda_{j}) } \left\{ \frac{1}{2|z_{-k}|} \sum_{i \in z_{-k}}(y_i-x_i^T\beta)^2 \right\},  \ \ \ \textnormal{if} \ \ \ \hat S^{(k)}(\lambda_{j}) \neq  \hat S^{(k)}(\lambda_{j-1}), \\
& \hat \beta^{(k)}_\textnormal{Lasso+OLS}(\lambda_{j-1}), \ \ \  \textnormal{otherwise};
\end{aligned}
\right.
\end{equation}

\item Compute the error in predicting the $k$th part of the data $PE^{(k)}$ :
\[ PE^{(k)}(\lambda_{j}) =  \frac{1}{|z_k|} \sum_{i \in z_k} \left( y_i - x_i^T \hat \beta^{(k)}_\textnormal{Lasso+OLS}(\lambda_{j}) \right)^2;\]
\end{itemize}

\STATE Compute cross validation error $CV(\lambda_{j})$, $j=1,...,J$:
\[ CV(\lambda_{j}) = \frac{1}{K} \sum_{k=1}^{K} PE^{(k)}(\lambda_{j}); \]
\STATE Compute the optimal $\lambda$ selected by CV;
\[ \lambda_{optimal} = \mathop {argmin}\limits_{\lambda_{j}:\ j=1,...,J} CV(\lambda_{j}); \]
\RETURN $\lambda_{optimal}$.
\end{algorithmic}
\end{algorithm*}

\begin{figure*}
\centerline{\includegraphics[width=.75\textwidth, height=0.4\textwidth]{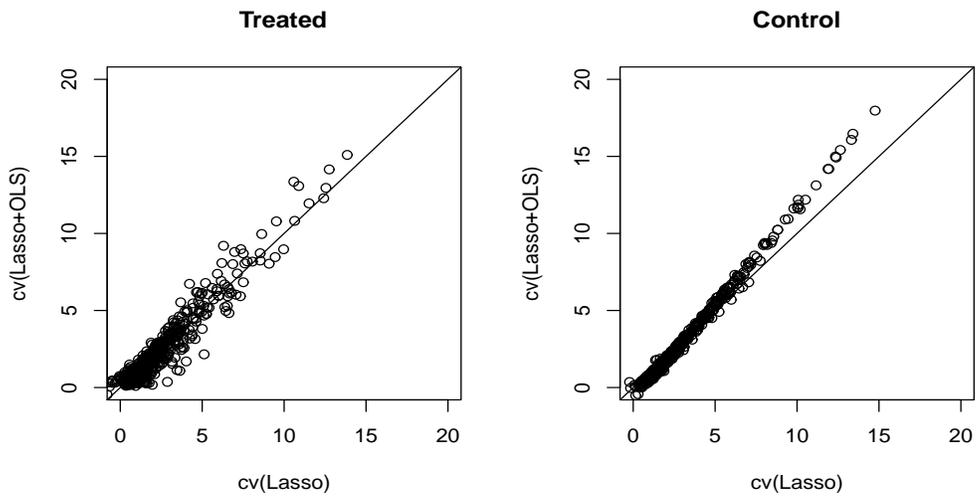}}
\vspace*{0.05in}
\caption{Adjustment (fitted) value comparison for cv(Lasso) and cv(Lasso+OLS).}\label{fig:fitvalue}
\end{figure*}

\begin{figure*}
\centerline{\includegraphics[width=.5\textwidth, height=0.5\textwidth]{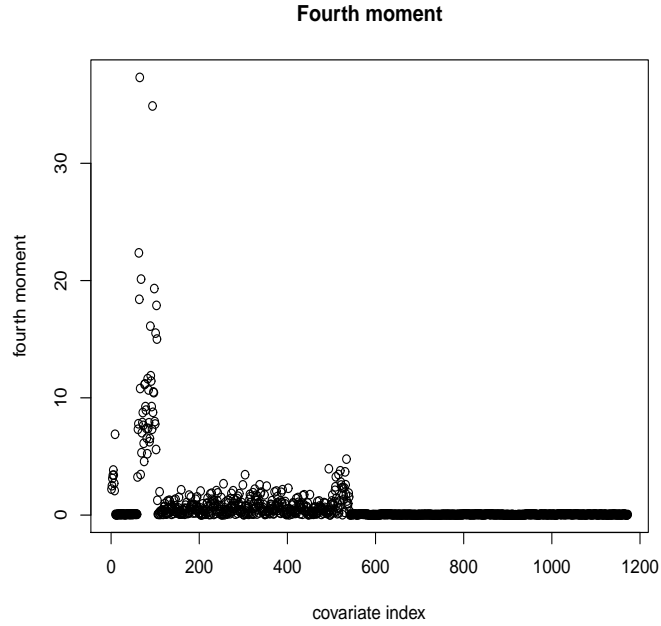}}
\vspace*{0.05in}
\caption{Fourth moment of each covariate. The covariates with the largest two fourth moments ($37.3$ and $34.9$ respectively) are quadratic term $interactnew^2$ and interaction term $IMscorerct:systemnew$ respectively. Neither of them are selected by the Lasso to do the adjustment. All the fourth moments of the main effects are less than $7$.}\label{fig:fourth-moment}
\end{figure*}

\begin{figure*}
\centerline{\includegraphics[width=.5\textwidth, height=0.5\textwidth]{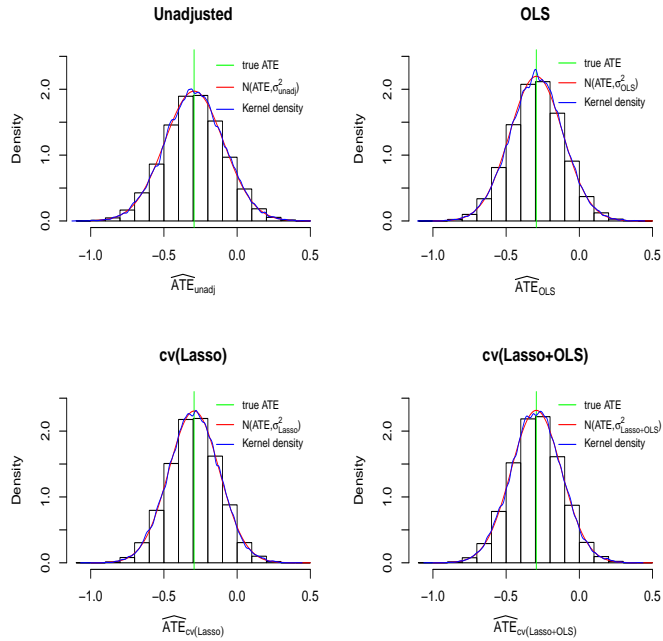}}
\vspace*{0.05in}
\caption{Histograms of ATE estimates. The green vertical lines are the true ATE; the red curves are the density of normal distribution; the blue curves are the kernel density estimate. The blue curves are very close to the red ones meaning that all the ATE estimates follow normal distribution.}\label{fig:density}
\end{figure*}


\begin{table*}
 \caption{\label{tab:mse} Bias, standard deviation (SD) and  root-mean square error $\sqrt{ \textnormal{MSE} }$ of ATE estimates}
\begin{tabular*}{\hsize}{@{\extracolsep{\fill}}llcccc}
 & & & & $(p,\rho)$ &  \\
  \cline{3-6}
 Statistic & Method  & (50,0) & (50,0.6) & (500,0) & (500,0.6) \\
  \hline  \\

  & & & $n_A=100$ & & \\ \hline
     & Unadjusted      & {\bf 0.003(0.004)} & {\bf 0.005(0.005)} & {\bf 0.002(0.003)} & {\bf 0.003(0.005)} \\
bias & OLS             & 0.014(0.005) & 0.013(0.006) & - & - \\
     & cv(Lasso)       & 0.007(0.004) & 0.014(0.005) & 0.006(0.004) & 0.005(0.004) \\
     & cv(Lasso+OLS)   & 0.011(0.004) & 0.013(0.005) & 0.009(0.004) & {\bf 0.003(0.004)} \\ \hline

     & Unadjusted      & 0.79(0.08) & 1.17(0.11) & 0.79(0.07) & 1.17(0.11) \\
SD   & OLS             & 0.72(0.07) & 0.96(0.09) & - & - \\
     & cv(Lasso)       & {\bf 0.62(0.06)} & {\bf 0.82(0.08)} & 0.67(0.06) & {\bf 0.84(0.08)} \\
     & cv(Lasso+OLS)   & 0.63(0.06) & {\bf 0.82(0.08)} & {\bf 0.65(0.06)} & {\bf 0.84(0.08)} \\ \hline
                            & Unadjusted      & 0.79(0.08) & 1.17(0.11) & 0.79(0.07) & 1.17(0.11) \\
$\sqrt{ \textnormal{MSE} }$ & OLS             & 0.72(0.07) & 0.97(0.09) & - & - \\
                            & cv(Lasso)       & {\bf 0.63(0.06)} & {\bf 0.82(0.08)} & 0.67(0.06) & 0.85(0.08) \\
                            & cv(Lasso+OLS)   & {\bf 0.63(0.06)} & {\bf 0.82(0.08)} & {\bf 0.65(0.06)} & {\bf 0.84(0.08)} \\ \hline
  \\
  & & & $n_A=125$ & & \\ \hline
     & Unadjusted      & 0.008(0.005) & 0.011(0.007) & {\bf 0.006(0.004)} & 0.01(0.007) \\
bias & OLS             & 0.008(0.004) & {\bf 0.005(0.005)} & - & - \\
     & cv(Lasso)       & {\bf 0.005(0.003)} & 0.012(0.005) & 0.007(0.004) & 0.004(0.004) \\
     & cv(Lasso+OLS)   & 0.012(0.004) & 0.012(0.005) & 0.011(0.004) & {\bf 0.003(0.003)} \\ \hline

     & Unadjusted      & 0.80(0.08) & 1.15(0.11) & 0.8(0.08) & 1.15(0.11) \\
SD   & OLS             & 0.69(0.06) & 0.90(0.09) & - & - \\
     & cv(Lasso)       & {\bf 0.62(0.06)} & {\bf 0.79(0.07)} & 0.67(0.06) & 0.82(0.08) \\
     & cv(Lasso+OLS)   & {\bf 0.62(0.06)} & {\bf 0.79(0.07)} & {\bf 0.65(0.06)} & {\bf 0.81(0.08)} \\ \hline
                            & Unadjusted      & 0.80(0.07) & 1.15(0.11) & 0.8(0.07) & 1.15(0.11) \\
$\sqrt{ \textnormal{MSE} }$ & OLS             & 0.69(0.07) & 0.90(0.09) & - & - \\
                            & cv(Lasso)       & {\bf 0.62(0.06)} & 0.80(0.08) & 0.67(0.06) & 0.82(0.08) \\
                            & cv(Lasso+OLS)   & {\bf 0.62(0.06)} & {\bf 0.79(0.07)} & {\bf 0.65(0.06)} & {\bf 0.81(0.08)} \\ \hline
  \\
  & & & $n_A=150$ & & \\ \hline
     & Unadjusted      & 0.004(0.004) & {\bf 0.000(0.005)} & {\bf 0.002(0.003)} & 0.005(0.005) \\
bias & OLS             & {\bf 0.002(0.003)} & 0.006(0.005) & - & - \\
     & cv(Lasso)       & 0.003(0.003) & 0.002(0.004) & 0.01(0.005) & 0.002(0.003) \\
     & cv(Lasso+OLS)   & 0.011(0.004) & 0.006(0.004) & 0.017(0.005) & {\bf 0.001(0.003)} \\ \hline

     & Unadjusted      & 0.85(0.08) & 1.19(0.11) & 0.85(0.08) & 1.19(0.11) \\
SD   & OLS             & 0.76(0.07) & 0.96(0.09) & - & - \\
     & cv(Lasso)       & {\bf 0.66(0.06)} & 0.82(0.08) & 0.72(0.07) & {\bf 0.84(0.08)} \\
     & cv(Lasso+OLS)   & 0.67(0.06) & {\bf 0.81(0.07)} & {\bf 0.71(0.07)} & {\bf 0.84(0.08)} \\ \hline

                            & Unadjusted      & 0.85(0.08) & 1.19(0.11) & 0.85(0.08) & 1.19(0.11) \\
$\sqrt{ \textnormal{MSE} }$ & OLS             & 0.76(0.07) & 0.96(0.09) & - & - \\
                            & cv(Lasso)       & {\bf 0.66(0.06)} & {\bf 0.82(0.08)} & 0.72(0.07) & {\bf 0.84(0.08)} \\
                            & cv(Lasso+OLS)   & 0.67(0.06) & {\bf 0.82(0.08)} & {\bf 0.71(0.07)} & {\bf 0.84(0.08)} \\
   \hline
\end{tabular*}
The numbers in parentheses are the corresponding standard errors estimated by using the bootstrap with $B=500$ resamplings of the ATE estimates.
\end{table*}

\begin{table*}
 \caption{\label{tab:modelsize} Mean number of selected covariates for treated and control group}
\begin{tabular*}{\hsize}{@{\extracolsep{\fill}}llcccc}
  \hline
 & & & & $(p,\rho)$ &  \\
  \cline{3-6}
 Group & Method  & (50,0) & (50,0.6) & (500,0) & (500,0.6) \\
  \hline  \\

   & & & $n_A=100$ & & \\ \hline
 treated & cv(Lasso)     & 16 & 13 & 22 & 22 \\
           & cv(Lasso+OLS) & 6 & 6 & 7 & 7 \\ \hline
 control   & cv(Lasso)     & 20 & 11 & 32 & 28 \\
           & cv(Lasso+OLS) & 8 & 6 & 7 & 7 \\ \hline
  \\
  & & & $n_A=125$ & & \\ \hline
 treated & cv(Lasso)     & 17 & 13 & 25 & 24 \\
           & cv(Lasso+OLS) & 7 & 6 & 6 & 6 \\ \hline
 control   & cv(Lasso)     & 19 & 11 & 32 & 27 \\
           & cv(Lasso+OLS) & 8 & 6 & 9 & 8 \\ \hline
  \\
  & & & $n_A=150$ & & \\ \hline
 treated & cv(Lasso)     & 18 & 13 & 29 & 26 \\
           & cv(Lasso+OLS) & 8 & 7 & 6 & 6 \\ \hline
 control   & cv(Lasso)     & 19 & 12 & 30 & 25 \\
           & cv(Lasso+OLS) & 8 & 6 & 11 & 8 \\
   \hline

\end{tabular*}
\end{table*}

\begin{table*}
 \caption{\label{tab:coverage} Coverage probability $(\%)$ and mean interval length (in parentheses) for $95\%$ confidence interval}
\begin{tabular*}{\hsize}{@{\extracolsep{\fill}}lcccc}
  \cline{2-5}
  & & $(p,\rho)$ & & \\
 Methods & (50,0) & (50,0.6) & (500,0) & (500,0.6) \\
  \hline  \\

  & & $n_A=100$ & &  \\
  \hline
  Unadjusted  & 97.3(3.54) & 95.8(4.79) & 97.3(3.54) & 95.8(4.79) \\
  OLS         & 92.2(2.55) & 90.0(3.19) & - & - \\
  cv(Lasso)       & 95.8(2.58) & 94.5(3.20) & 94.3(2.61) & 92.4(3.07) \\
  cv(Lasso+OLS)   & 95.6(2.57) & 94.4(3.17) & 94.8(2.60) & 93.0(3.11) \\ \hline
  \\
  & & $n_A=125$ & & \\
  \hline
  Unadjusted  & 97.4(3.56) & 96.0(4.74) & 97.3(3.56) & 95.9(4.74) \\
  OLS         & 93.3(2.54) & 91.6(3.14) & - & - \\
  cv(Lasso)       & 96.0(2.56) & 95.0(3.15) & 94.1(2.59) & 92.9(3.02) \\
  cv(Lasso+OLS)   & 95.7(2.55) & 94.9(3.12) & 94.4(2.58) & 93.6(3.06) \\ \hline
  \\
  & & $n_A=150$ & &  \\
  \hline
  Unadjusted  & 97.1(3.72) & 95.8(4.88) & 97.1(3.72) & 95.8(4.88) \\
  OLS         & 91.4(2.64) & 90.4(3.21) & - & - \\
  cv(Lasso)       & 95.4(2.66) & 94.9(3.23) & 92.9(2.68) & 92.6(3.08) \\
  cv(Lasso+OLS)   & 94.7(2.63) & 94.8(3.19) & 92.0(2.63) & 93.1(3.11) \\
   \hline
\end{tabular*}
The numbers in parentheses are the corresponding mean interval lengths.
\end{table*}

\end{document}